\documentclass[10pt]{amsart}
\usepackage{amsmath}
\usepackage{latexsym,amsfonts,amssymb,mathrsfs}
\usepackage{mathdots}
\usepackage{color}
\usepackage[all,2cell,color]{xy}
\usepackage{enumitem} 
\usepackage{bbm}
\usepackage{tikz}
\usetikzlibrary{cd}
\usetikzlibrary{calc}
\usetikzlibrary{math}
\usetikzlibrary{arrows}
\usetikzlibrary{fit}
\usetikzlibrary{positioning}
\usetikzlibrary{decorations.pathmorphing, arrows.meta}
\tikzset{Rightarrow/.style={double equal sign distance,>={Implies},->},
triple/.style={-,preaction={draw,Rightarrow}},
quadruple/.style={preaction={draw,Rightarrow,shorten >=0pt},shorten >=1pt,-,double,double
distance=0.2pt}}

\usepackage[T1]{fontenc}
 \usepackage[utf8]{inputenc}
 \usepackage{array}
 \usepackage{blkarray} 
 \usepackage{mathtools} 
 \usepackage{longtable}
 \usepackage{multicol} 
 \usepackage{hyperref}
\usepackage[capitalise]{cleveref}

\usetikzlibrary{positioning} 
\usetikzlibrary{decorations.pathreplacing}
\usetikzlibrary{arrows, decorations.markings} 

\newcolumntype{D}{>{\hfil$}p{1.2cm}<{$\hfil}}

\tikzset{%
    symbol/.style={%
        draw=none,
        every to/.append style={%
            edge node={node [sloped, allow upside down, auto=false]{$#1$}}}
    }
}

\tikzset{%
scalearrow/.style n args={3}{
  decoration={
    markings,
    mark=at position (1-#1)/2*\pgfdecoratedpathlength
      with {\coordinate (#2);},
    mark=at position (1+#1)/2*\pgfdecoratedpathlength
      with {\coordinate (#3);},
    },
  postaction=decorate,
  }
}

\theoremstyle{plain}   
\newtheorem{thm}{Theorem}[section] 
\makeatletter\let\c@thm\c@thm\makeatother

\makeatletter\let\c@cor\c@thm\makeatother
\newtheorem{lem}{Lemma}[section]
\makeatletter\let\c@lem\c@thm\makeatother
\newtheorem{prop}{Proposition}[section]
\makeatletter\let\c@prop\c@thm\makeatother

\makeatletter\let\c@claim\c@thm\makeatother

\makeatletter\let\c@conjecture\c@thm\makeatother

\makeatletter\let\c@wconjecture\c@thm\makeatother

\newtheorem{thmalph}{Theorem}

\theoremstyle{definition}

\newtheorem{defn}{Definition}[section]
\makeatletter\let\c@defn\c@thm\makeatother
\newtheorem{const}{Construction}[section]
\makeatletter\let\c@const\c@thm\makeatother
\newtheorem{notn}{Notation}[section]
\makeatletter\let\c@notn\c@thm\makeatother

\makeatletter\let\c@convention\c@thm\makeatother

\makeatletter\let\c@convention\c@thm\makeatother

\theoremstyle{remark}

\newtheorem{rmk}{Remark}[section]
\makeatletter\let\c@rmk\c@thm\makeatother
\newtheorem{ex}{Example}[section]
\makeatletter\let\c@ex\c@thm\makeatother

\makeatletter\let\c@observation\c@thm\makeatother

\makeatletter\let\c@warning\c@thm\makeatother

\makeatletter\let\c@digression\c@thm\makeatother

\makeatletter\let\c@answ\c@thm\makeatother

\makeatletter\let\c@answ\c@thm\makeatother

\makeatletter\let\c@aside\c@thm\makeatother

\makeatletter
\let\c@equation\c@thm
\numberwithin{equation}{section}
\makeatother

\crefname{lem}{Lemma}{Lemmas}
\crefname{thm}{Theorem}{Theorems}
\crefname{defn}{Definition}{Definitions}
\crefname{notn}{Notation}{Notations}
\crefname{const}{Construction}{Constructions}
\crefname{prop}{Proposition}{Propositions}
\crefname{rmk}{Remark}{Remarks}
\crefname{cor}{Corollary}{Corollaries}
\crefname{equation}{Display}{Displays}
\crefname{ex}{Example}{Examples}
\crefname{thmalph}{Theorem}{Theorems}
\crefname{answ}{Answer}{Answers}
\crefname{question}{Question}{Questions}

\newcommand{\cA}{\mathcal{A}}

\newcommand{\cC}{\mathcal{C}}
\newcommand{\cD}{\mathcal{D}}

\newcommand{\cM}{\mathcal{M}}

\newcommand{\cO}{\mathcal{O}}

\newcommand{\cS}{\mathcal{S}}

\newcommand{\bT}{\Theta_n}

\newcommand{\cat}{\cC\!\mathit{at}}
\newcommand{\set}{\cS\!\mathit{et}}
\newcommand{\sset}{\mathit{s}\set}

\newcommand{\adch}{\mathit{ad}\mathcal C\mathit{h}}

\newcommand{\psh}[1]{\set^{#1^{\op}}}
\newcommand{\spsh}[1]{\sset^{#1^{\op}}}
\newcommand{\msset}{m\sset}

 \newcommand{\omegacat}{\omega\cat}

 \newcommand{\comp}{\ast}

\newcommand{\inttrunc}[1]{\tau_{\leq n}^{\text{i}}}

\DeclareMathOperator{\id}{id}

\newcommand{\aamalg}[1]{\underset{#1}{{\amalg}}} 
\DeclareMathOperator{\op}{op}

\newcommand{\NRS}{N^{\operatorname{RS}}}

\newcommand{\Gray}{\otimes}
\newcommand{\chsusp}{\Sigma}

\newcommand{\tabld}[2]{\begin{pmatrix}#1^-_0 &\dots &#1^-_{#2-1}
  &#1^-_{#2}\cr\noalign{\vskip 3pt} #1^+_0 &\dots &#1^+_{#2-1}
  &#1^+_{#2}\end{pmatrix}}

\newcommand{\suspindex}{r}

\AtBeginDocument{%
   \def\MR#1{}
}

\usepackage[margin=1.55in]{geometry}

\keywords{$(\infty,n)$-categories, $n$-categories, complicial sets, complete Segal $\Theta_n$-spaces, suspension of $n$-categories, complicial nerve}
\subjclass[2020]{18N65; 55U35; 18N50; 55U10}

\begin{document}

\author{Viktoriya Ozornova}
\address{Max Planck Institute for Mathematics, Bonn, Germany}
\email{viktoriya.ozornova@mpim-bonn.mpg.de}

\author{Martina Rovelli}
\address{Department of Mathematics and Statistics, 
University of Massachusetts, MA 01003-9305
Amherst, USA
}
\email{rovelli@math.umass.edu} 

\title[Globular and complicial approaches to $(\infty,n)$-categories]{A Quillen adjunction between\\ globular and complicial approaches to $(\infty,n)$-categories}

\maketitle

\begin{abstract}
We prove the compatibility between the suspension construction and the complicial nerve of $\omega$-categories. As a motivating application, we produce a Quillen pair between the models of $(\infty,n)$-categories given by Rezk's complete Segal $\Theta_n$-spaces and Verity's $n$-complicial sets.
    \end{abstract}

\tableofcontents

\section*{Introduction}

It is by now known that many mathematical phenomena of interest can only be properly formalized using the language of $(\infty,n)$-categories. Several mathematical objects have been identified to implement the notion of an $(\infty,n)$-category, each with its own advantages and disadvantages. Amongst those, there are Verity's $n$-complicial sets \cite{VerityComplicialI,VeritySlides,RiehlNotes,OR,RiehlVerityBook} and Rezk's complete Segal $\Theta_n$-spaces \cite{rezkTheta}.

The homotopy theories of $n$-complicial sets and complete Segal $\Theta_n$-spaces are only known to be equivalent for $n\le2$, and this paper reports progress towards establishing the equivalence of these homotopy theories for general $n$, which was conjectured more than three decades ago (see e.g.~\cite{StreetOrientedSimplexes,VeritySlides,BarwickSchommerPries}).

\begin{thmalph}
   There is an adjunction of $\infty$-categories
between the $\infty$-category of complete Segal $\Theta_n$-spaces and the $\infty$-category of $n$-complicial sets.
\end{thmalph}
More precisely, we achieve this by constructing an adjunction between the model categories $\spsh{\Theta_n}_{(\infty,n)}$  and $\msset_{(\infty,n)}$, which we show in \cref{LQuillenPair} to be a Quillen pair (and which is conjecturally a Quillen equivalence). This results (partially) generalizes joint work of the authors with Bergner \cite{BOR}.

In order to prove \cref{LQuillenPair}, the crucial ingredient is to understand how the two-point suspension interacts with the complicial nerve of certain $n$-categories. We prove the following as \cref{ThmB}.
        \begin{thmalph}
If $\cC$ admits an algebraic model in an appropriate sense, then $\NRS\Sigma\cC$ is equivalent to $\Sigma\NRS\cC$ in the model structure for $n$-complicial sets. 
    \end{thmalph}
Here, the precise condition on $\cC$ requires it to be obtained from an (augmented directed) chain complex via Steiner's functor $\nu\colon\adch\to\omega\cat$ (see \cite{SteinerEmbedding,AraMaltsiniotisJoin}), the functor $\NRS\colon n\cat\to\msset$ is the \emph{Roberts--Street nerve} and the functors $\Sigma\colon n\cat\to(n+1)\cat$ and $\Sigma\colon\msset_{(\infty,n)}\to\msset_{(\infty,n+1)}$ implement the two-point suspension construction in a strict and weak context.
The theorem is also used in work by Loubaton \cite{loubaton2}, who gives a criterion to identify self-equivalences on the $\infty$-category of $n$-complicial sets.

\addtocontents{toc}{\protect\setcounter{tocdepth}{1}}
\subsection*{Acknowledgements}
It is hard to overestimate the role of Andrea Gagna for this paper, who has taught the authors the language of Steiner's theory of augmented directed chain complexes, without which the current result would have been out of our reach. This work was completed while the authors
visited the 
 Instituto de Matemáticas de UNAM in Cuernavaca for the program \emph{Higher categories -- Part 2}, supported by the National Science Foundation under Grant No. DMS-1928930. The second author is grateful for support from the National Science
Foundation under Grant No. DMS-2203915.

\tableofcontents

\section{Steiner's augmented directed chain complexes}

We recall the basic definitions around Steiner's augmented directed chain complexes, as well as some constructions based on augmented directed chain complexes: the suspension, tensor product, and the total dual, as well as the main properties that we use later in the paper and relevant examples. Most of the material is drawn from \cite{SteinerEmbedding} (see also \cite{AraMaltsiniotisJoin}).

\subsection{Augmented directed chain complexes}
\label{sec:Steiner}

By a \emph{chain complex} $C$ we will always mean an
$\mathbb{N}$-graded chain complex of abelian
groups with homological indexing, that is,
a family $(C_q)_{q\geq 0}$ of abelian groups, together with maps $\partial_{q} \colon C_{q+1} \to C_{q}$ satisfying $\partial_q \partial_{q+1}=0$. We also assume that, whenever occurring, $C_{-1}=0$, and $\partial_{0}=0$.

Given chain complexes $C$ and $\overline C$, a \emph{chain map} or \emph{morphism of chain complexes} $\phi\colon C\to \overline C$ consists of a family of homomorphisms $(\phi_q\colon C_q\to \overline C_q)_{q\geq 0}$ that commutes with the differentials in the sense that $\overline\partial_{q} \phi_{q+1}=\phi_q\partial_{q}$
for every $q \geq 0$.

An \emph{augmented chain complex} is a pair $(C,\varepsilon)$
of a chain complex $C$ and an augmentation, namely a map $\varepsilon\colon C_0\to\mathbb Z$ such that $\varepsilon\partial_0=0$.

An \emph{augmented chain map} $\phi \colon (C,\varepsilon) \to (\overline C, \overline\varepsilon)$ between augmented chain complexes $(C,\varepsilon)$ and $(\overline C,\overline\varepsilon)$ consists of a chain map $\phi\colon C\to\overline C$ that is moreover compatible with the augmentations, namely such that $\overline\varepsilon \phi_0 = \varepsilon.$

We recall the enhancement of augmented chain complexes developed by Steiner \cite[\textsection2]{SteinerEmbedding}.

\begin{defn}[{\cite[Def.~2.2]{SteinerEmbedding}}]
 An \emph{augmented directed complex}
 is a triple
 $(C, C^+, \varepsilon)$ where $(C, \varepsilon)$ is
 an augmented chain complex and $C^+ =  (C^+_q)_{q\ge0}$ is a collection of commutative monoids,
 where $C^+_q$ is a submonoid of $C_q$ called the
 \emph{positivity submonoid} of $C_q$.
 
 A \emph{morphism of augmented directed chain complexes}, or an \emph{augmented directed chain map} $\phi \colon (C, C^+, \varepsilon) \to (\overline C,  \overline C^+, \overline\varepsilon)$ between augmented directed chain complexes $(C, C^+, \varepsilon)$ and $( \overline C,  \overline C^+, \overline\varepsilon)$
 is an augmented chain map $\phi \colon (C, \varepsilon) \to (\overline C, \overline{\varepsilon})$ that moreover preserves the positivity submonoids,
    namely such that
    \[\phi_q(C^+_q) \subseteq  \overline C^+_q\]
 for all $q\ge0$.
\end{defn}

We denote by $\adch$ the category of augmented directed chain complexes and maps of chain complexes that preserve the augmentation and the positivity submonoids.

\begin{rmk}
\label{ColimitsADCH}
The category $\adch$ is cocomplete, colimits are computed degreewise, and epimorphisms are detected pointwise in the category $\cA b$ of abelian groups and the category $c\cM on$ of commutative monoids, which are both cocomplete. That is, the forgetful functor
\[\adch\to\prod_{q\ge0}(\cA b\times c\cM on)\]
given by $C\mapsto(C_q,C_q^+)_{q\ge0}$ creates colimits (and in particular epimorphisms).
\end{rmk}

\begin{rmk}
\label{ColimitsADCH2}
Consider the following left adjoint functors.
\begin{enumerate}[leftmargin=*]
    \item The free abelian group functor on a set and the free commutative monoid functor on a set,
\[\mathbb Z[-]\colon\set\to\cA b\text{ and }\mathbb N[-]\colon\set\to c\cM on,\] given by $X\mapsto\mathbb Z[X]$ and $X\mapsto\mathbb N[X]$. The right adjoint functors are the forgetful functors.
\item The free abelian group functor on a pointed set and the free commutative monoid functor on a pointed set,
\[\mathbb Z[-]\colon\set_*\to\cA b\text{ and }\mathbb N[-]\colon\set_*\to c\cM on,\] given by $(X,x_0)\mapsto\mathbb Z[X\setminus\{x_0\}]$ and $(X,x_0)\mapsto\mathbb N[X\setminus \{x_0\}]$. The right adjoint functors are the forgetful functors that retain the identity as a base point.
\item The functor that freely adds a base point to a set, \[(-)_+\colon\set\to\set_*,\]
given by $X\mapsto (X\amalg\{*\},*)$. The right adjoint functor is the functor that forgets the base point.
\end{enumerate}
Being left adjoint functors, they all preserve colimits (and in particular epimorphisms).
\end{rmk}

\begin{notn} Let $m\ge-1$ and $q\ge -1$. We denote
\begin{itemize}[leftmargin=*]
    \item by $\Delta[m]_q=\cat([q],[m])$ the set of $q$-simplices of the standard simplex\footnote{We follow the convention that $[-1]$ is the empty category, and $\Delta[-1]$ is the initial simplicial set, which is levelwise empty.} $\Delta[m]$. A generic $q$-simplex in $\Delta[m]$ is of the form
\[
[\mathbf a]=[a_0,a_1,\dots,a_q]
\]
with $0\leq a_0 \le a_1 \le \ldots \le a_q\leq m$. We say that $q$ is the \emph{length} $|\mathbf{a}|$ of $[\mathbf{a}]$.
\item by $B[m]_q\subseteq\Delta[m]_q $ the set of non-degenerate $q$-simplices of $\Delta[m]$, namely those simplices for which $0\leq a_0 < a_1 < \ldots < a_q\leq m$.
\item by $O[m]_q=\mathbb Z[B[m]_q]\cong\mathbb Z^{{[m]}\choose{[q]}}$ the abelian group freely generated by non-degenerate $q$-simplices of $\Delta[m]$.
The generic element of $O[m]_q$ is a formal sum
\[
        c = \sum_{[\mathbf a] \in B[m]_q} c_{[\mathbf a]} \cdot [\mathbf a]
    \]
    where $c_{[\mathbf a]} \in \mathbb{Z}$.
\item by $O[m]_q^+=\mathbb N[B[m]_q]\cong\mathbb N^{{[m]\choose[q]}}$ the abelian monoid freely generated by non-degenerate $q$-simplices of $\Delta[m]$. The generic element of $O[m]_q^+$ is one for which $c_{[\mathbf a]} \in \mathbb{N}$.
\end{itemize}
There are canonical inclusions $\Delta[m]_q\supseteq B[m]_q\subseteq O[m]_q^+\subseteq O[m]_q$.
\end{notn}

The augmented directed chain complex $O[m]$ is the algebraic model of the $m$-oriental $\cO[m]$, in a sense that will be made precise in \cref{orientalandnu}.

\begin{ex}[{\cite[Ex.~3.8]{SteinerEmbedding}}]
\label{AlgebraicOrientals}
For $m\ge-1$, we consider the augmented directed chain complex
$O[m]$ with the following structure.
\begin{itemize}[leftmargin=*]
    \item For $q\ge0$ the abelian group of $q$-chains is given by $O[m]_q$.
    \item For $q\ge0$ the commutative monoid of positive $q$-chains is given by $O[m]_q^+$.
    \item For $q\ge-1$ the differential $\partial_q\colon O[m]_{q+1}\to O[m]_{q}$ is given by
\[
\begin{array}{llllll}
     \partial_q[\mathbf a]&=\partial_q[a_0,\dots,a_q,a_{q+1}]  \\
     & =\sum\limits_{i=0}^{q+1}(-1)^i\cdot[a_0,\dots,\widehat{a_i},\dots,a_{q+1}]\in O[m]_{q}\\
\end{array}
\]
\item The augmentation map $\varepsilon\colon O[m]_0\to\mathbb Z$ is given by
\[\varepsilon[a]=1\in\mathbb Z.\]
\end{itemize}
\end{ex}

Later in the paper, we will make use of the following dual construction for an augmented directed chain complex.

\begin{defn}[{\cite[\textsection2.18]{AraMaltsiniotisJoin}}]
\label{defn:op_ADC}
Let $C$ be an augmented directed complex.
    The \emph{total dual} $C^{\circ}$ of $C$
    is the augmented directed complex
    with the following structure
   \begin{itemize}[leftmargin=*]
        \item For $q\ge1$ the abelian group of $q$-chains is given by $C^{\circ}_q = C_q$;
        \item For $q\ge1$ the commutative monoid of positive chains is given by $(C^{\circ})_q^+ = C^+_q$;
        \item For $q\ge1$, the differential $\partial^{C^\circ}_q\colon C^\circ_q\to C^\circ_{q-1}$ is given by $\partial^{C^\circ}_q(c)=-\partial^{C}_q(c)$.
        \item The augmentation $\varepsilon^{C^\circ}\colon C^\circ_0\to \mathbb Z$ is given by $\varepsilon^{C^\circ}(a)=\varepsilon^{C}(a)$.
    \end{itemize}
    This construction defines an involution $(-)^{\circ} \colon \adch \to \adch$.
\end{defn}

\subsection{Suspension of augmented directed chain complexes}

We define a two-point suspension for augmented directed chain complexes. This is the construction that Steiner denotes $V(1,C)$ in \cite[\textsection5]{SteinerSimpleOmega}\footnote{This is different from the one-point suspension considered by Ara--Maltsiniotis in \cite[\textsection6.3]{AraMaltsiniotisJoin}.}.

\begin{defn}
Let $C$ be an augmented directed chain complex. The \emph{suspension}  of $C$ is the augmented directed chain complex $\Sigma C$ with the following structure:
\begin{itemize}[leftmargin=*]
    \item For $q\ge0$, the abelian group $(\Sigma C)_q$ of $q$-chains is given by
\[(\Sigma C)_q=\left\{\begin{array}{cccc} \mathbb{Z}[\bot,\top] &\mbox{ if }q=0,\\
    C_{q-1} &\mbox{ if }q\ge1.
    \end{array}\right.\]
    \item For $q\ge0$, the commutative monoid $(\Sigma C)_q^+$ of positive $q$-chains is given by
    \[(\Sigma C)_q^+=\left\{\begin{array}{cccc} \mathbb{N}[\bot,\top] &\mbox{ if }q=0,\\
    C_{q-1}^+ &\mbox{ if }q\ge1.
    \end{array}\right.\]

        \item For $q\ge0$, the differential $\partial_q\colon(\Sigma C)_q\to(\Sigma C)_{q-1}$ is given by
     \[\partial_q^{\Sigma C}(a):=\left\{
        \begin{array}{llll}
        \varepsilon^C(c)\cdot(\top-\bot)=-\varepsilon^C(c)\cdot \bot+\varepsilon^C(c)\cdot \top&\mbox{ if }q=0,\\
        \partial^C_{q-1}(c)&\mbox{ if }q\ge1.
        \end{array}\right.\]
\item The augmentation $\varepsilon^{\Sigma C}\colon(\Sigma C)_0\to\mathbb Z$ is given by
\[ \varepsilon^{\Sigma C}\bot=1=\varepsilon^{\Sigma C}\top.\]
\end{itemize}

\end{defn}

The augmented directed chain complex $\Sigma C$ comes with a map $\Sigma O[-1]\to\Sigma C$ so it can be naturally regarded as an object of $\prescript{{\Sigma O[-1]/}}{}{\adch}$. The following is a consequence of \cite[Theorem 5.6]{SteinerSimpleOmega}.

\begin{prop}
\label{suspensionFF}
The suspension functor $\Sigma\colon\adch\to\prescript{{\Sigma O[-1]/}}{}{\adch}$ is fully faithful. 
\end{prop}

\subsection{Tensor product of augmented directed chain complexes}

We consider the tensor product of abelian groups $\otimes\colon\cA b\times\cA b\to\cA b$, as well as the (less known) tensor product of commutative monoids $\otimes\colon c\cM on\times c\cM on\to c\cM on$. See e.g.~\cite[Chapter 16]{GolanSemiringsAppl}
for more details on this construction. We will mostly use instances of the tensor product of \emph{free} abelian groups and \emph{free} commutative monoids, which is described by the following remark.

\begin{rmk}
Recall the functors from \cref{ColimitsADCH2}.
\begin{enumerate}[leftmargin=*]
    \item The free abelian group functor $\mathbb Z[-]\colon(\set,\times)\to(\cA b,\otimes)$ and the free commutative monoid functor $\mathbb N[-]\colon(\set,\times)\to (c\cM on,\otimes)$ is strong monoidal, namely, there are natural bijections
\[\mathbb Z[X]\otimes\mathbb Z[Y]\cong \mathbb Z[X\times Y]\quad\text{ and }\quad\mathbb N[X]\otimes\mathbb N[Y]\cong \mathbb N[X\times Y],\]
for any $X$ and $Y$ sets.
In particular, the tensor product of free abelian groups, resp.~commutative monoids, is a free abelian group, resp.~free commutative monoid.
It follows that also the free abelian group functor $\mathbb Z[-]\colon(\set_*,\times)\to(\cA b,\otimes)$ and the free commutative monoid functor $\mathbb N[-]\colon(\set,\times)\to (c\cM on,\otimes)$ is strong monoidal
\item The free abelian group functor $\mathbb Z[-]\colon(\set,\amalg)\to(\cA b,\oplus)$ and the free commutative monoid functor $\mathbb N[-]\colon(\set,\amalg)\to (c\cM on,\oplus)$ is strong monoidal. 
\end{enumerate}
\end{rmk}

\begin{defn}[{\cite[Example 3.10]{SteinerEmbedding}}]
\label{tensoradch}
Let $C$ and $ D$ be augmented directed chain complexes. The \emph{tensor product} of $C$ and $ D$ is the augmented directed chain complex $C\otimes  D$ with the following structure:
\begin{itemize}[leftmargin=*]
    \item For $q\ge0$, the abelian group $(C\otimes  D)_q$ of $q$-chains is given by
\[(C\otimes D)_q=\bigoplus_{k+\ell=q}C_k\otimes D_{\ell}.\]
    \item For $q\ge0$, the commutative monoid $(C\otimes  D)_q^+$ of positive $q$-chains is given by
    \[(C\otimes  D)_q^+=\bigoplus_{k+\ell=q}C_k^+\otimes D_{\ell}^+\]
    \item For $q\ge0$, the differential $\partial_q^{C\otimes  D}\colon(C\otimes  D)_q\to(C\otimes  D)_{q-1}$ is given by
        \[\partial_q^{C\otimes  D}(c\otimes d):=\partial^C c\otimes D+(-1)^{|c|} c\otimes\partial^{ D} d
       \]
\item The augmentation $\varepsilon^{C\otimes  D}\colon(C\otimes  D)_0\cong C_0\otimes D_0\to\mathbb Z$ is given by
\[ \varepsilon^{C\otimes  D}(c\otimes  d)=\varepsilon ^Cc\cdot\varepsilon^D d.\]
\end{itemize}
\end{defn}

The construction defines a functor $\otimes\colon\adch\times\adch\to \adch$.

We now unpack tensor product of orientals.

\begin{ex}
Let $k,\ell\ge0$.
\begin{itemize}[leftmargin=*]
    \item For $q\ge0$, the abelian group $(O[k]\otimes O[\ell]^\circ)_q$ of $q$-chains is given by  
    \[(O[k]\otimes O[\ell]^\circ)_q=\bigoplus_{i=0}^{q}O[k]_i\otimes O[\ell]_{q-i}\cong\bigoplus\limits_{i=0}^q\mathbb Z[B[k]_i]\otimes\mathbb Z[B[\ell]_{q-i}]\cong\bigoplus\limits_{i=0}^q\mathbb Z^{[k]\choose[i]}\otimes\mathbb Z^{[\ell]\choose[q-i]}\]
    \item For $q\ge0$, the commutative monoid $(O[k]\otimes O[\ell]^\circ)_q^+$ of positive $q$-chains is given by
    \[(O[k]\otimes O[\ell]^\circ)_q^+=\bigoplus_{i=0}^{q}O[k]^+_i\otimes O[\ell]^+_{q-i}\cong\bigoplus\limits_{i=0}^q\mathbb N[B[k]_i]\otimes\mathbb N[B[\ell]_{q-i}]\cong\bigoplus\limits_{i=0}^q\mathbb N^{[k]\choose[i]}\otimes\mathbb N^{[\ell]\choose[q-i]}\]

        \item For $q>0$, the differential
    $\partial_q^{O[k]\otimes O[\ell]^\circ}\colon(O[k]\otimes O[\ell]^\circ)_{q+1}\to(O[k]\otimes O[\ell]^\circ)_{q}$ is given by
        \[\partial_q^{O[k]\otimes O[\ell]^\circ}([\textbf a]\otimes [\textbf b]):=\partial^{O[k]} [\textbf a]\otimes [\textbf b]+(-1)^{|\textbf a|} [\textbf a]\otimes\partial^{ O[\ell]^{\circ}} [\textbf b]
       \]
\item The augmentation $\varepsilon^{O[k]\otimes O[\ell]^\circ}\colon(O[k]\otimes O[\ell]^\circ)_0\cong O[k]_0\otimes O[\ell]^\circ_0\to\mathbb Z$ is given by
\[ \varepsilon^{O[k]\otimes  O[\ell]^\circ}([a]\otimes  [b])\coloneqq\varepsilon ^{O[k]}([a])\cdot\varepsilon^{ O[\ell]^\circ}([b]) =1\cdot 1=1\]
\end{itemize}
\end{ex}

Recall the functor $(-)_+\colon\set\to\set_*$, which is the left adjoint to the forgetful functor.

\begin{rmk}
\label{descriptionphi}
For $q\ge 0$, there is a canonical map of pointed sets
\[
{[k+1+\ell]\choose [q]}_+ \to \bigvee\limits_{r=0}^q \left({[k]\choose[r]}\times{[\ell]\choose[q-r]}\right)_+
\]
that splits a subset of $[k+1+\ell]$ into its $[k]$-part and $[\ell]$-part, using the base point whenever any of them is empty. This induces a map of abelian groups 
\[\phi_q\colon O[k+1+\ell]_q\cong\mathbb Z^{[k+1+\ell]\choose [q]}\to\bigoplus\limits_{r=0}^q\mathbb Z^{[k]\choose[r]}\otimes \mathbb{Z}^{[\ell]\choose[q-r]}\cong \Sigma(O[k]\otimes O[\ell]^\circ)_q\]
and commutative monoids
\[\phi_q\colon O[k+1+\ell]^+_q\cong\mathbb N^{{[k+1+\ell]\choose [q]}}\to \bigoplus\limits_{r=0}^q\mathbb{N}^{[k]\choose[r]}\otimes \mathbb{N}^{[\ell]\choose[q-r]}\cong \Sigma(O[k]\otimes O[\ell]^\circ)_q^+.\]
Explicitly, $\phi_0\colon O[k+1+\ell]_0\to \chsusp\left(O[k]\otimes O[\ell])^{\circ}\right)_0$ is given by
\[\phi_0([a'']):=\left\{
\begin{array}{llll}
     \bot& \text{if }0\leq a'' \leq k, \\
   \top  & \text{if }k+1 \leq a'' \leq k+1+\ell.
\end{array}
\right.\]
For $q>0$, the map $\phi_q\colon O[k+1+\ell]_q\to \chsusp\left(O[k]\otimes O[\ell]^{\circ}\right)_q$ is given by
\[\phi_q([\mathbf{a}, \mathbf{a'}]):=
\left\{
\begin{array}{llll}
[\mathbf{a}]\otimes (s^0)^{k+1}[\mathbf{a'}], &\text{ with }
\mathbf{a}\subseteq [0,k], \mathbf{a'}\subseteq [k+1,m], |\mathbf a|\ge0, |\mathbf a'|\ge0, \\
   0  & \text{else.}
\end{array}
\right.\]
\end{rmk}

\begin{prop}
\label{orientaltosuspension}
Let $k,\ell\ge0$. There is a map of augmented directed chain complexes
\begin{equation}
   \label{orientaltosuspension1}
   \phi\colon O[k+1+\ell]\to \chsusp\left(O[k]\otimes O[\ell]^{\circ}\right).\end{equation}
   given degreewise by the maps $\phi_q$ described in \cref{descriptionphi}.
\end{prop}

\begin{proof}
Given $\mathbf{a}\subseteq [0,k]$, and $\mathbf{a'}\subseteq [k+1,n]$, with $|\mathbf a|\ge0$ and $|\mathbf a'|\ge0$, which is the only case of interest, we obtain
\[
\begin{array}{lllll}
\partial\phi([\mathbf{a}, \mathbf{a'}])
&=&&\partial\left([\mathbf{a}]\otimes (s^0)^{k+1}[\mathbf{a'}]\right)\\
&=&&\partial[\mathbf{a}]\otimes (s^0)^{k+1}[\mathbf{a'}] + (-1)^{|\mathbf a|+1}[\mathbf{a}]\otimes \partial (s^0)^{k+1}[\mathbf{a'}]
  \\
 &=&&\phi(\partial([\mathbf{a}, \mathbf{a'}]))\\
\end{array}
\]
as desired.
\end{proof}

For $k,\ell\geq 0$,
the maps from \cref{orientaltosuspension}, together with other canonical maps, can be used to build a commutative diagram
    \begin{equation}
\label{sigmaoforientaltensor}
\begin{tikzcd}
 O[k]\oplus  O[\ell]^{\circ} \arrow[r]\arrow[d] &  O[k+1+\ell]\arrow[d]\\
 O[0]\oplus  O[0]^\circ\arrow[r] & \Sigma( O[k]\Gray  O[\ell]^{\circ}).
\end{tikzcd}
    \end{equation}

\begin{prop}
 \label{quotientoforiental} Let $k,\ell\geq 0$.
 The diagram \eqref{sigmaoforientaltensor} induces a natural isomorphism of augmented directed chain complexes
 \[\Sigma( O[k]\Gray  O[\ell]^{\circ})\cong(O[0]\oplus  O[0]^\circ)\aamalg{O[k]\oplus  O[\ell]^\circ}O[k+1+\ell].\]
\end{prop}

\begin{proof}
There are pushout squares of abelian groups and of commutative monoids:
 \[
\begin{tikzcd}
 \mathbb Z^{[k]}\oplus \mathbb Z^{[\ell]}\arrow[r,"\cong"]\arrow[d] &\mathbb Z^{[k+1+\ell]}\arrow[d]\\
\mathbb Z\oplus\mathbb Z\arrow[r,"\cong" swap] &\mathbb Z[\bot,\top]
\end{tikzcd}
\quad\text{ and }\quad
\begin{tikzcd}
 \mathbb N^{[k]}\oplus \mathbb N^{[\ell]}\arrow[r,"\cong"]\arrow[d] &\mathbb N^{[k+1+\ell]}\arrow[d]\\
\mathbb N\oplus\mathbb N\arrow[r,"\cong" swap] &\mathbb N[\bot,\top].
\end{tikzcd}
\]
Here, the left vertical maps are the sums of the canonical map that folds the first $k+1$ copies of $\mathbb Z$ and the one that folds the last $\ell+1$ copies of $\mathbb Z$.
This means that \eqref{sigmaoforientaltensor} induces a pushout of abelian groups (resp., commutative monoids):
    \begin{equation} \label{sigmaoforientaltensor0}
\begin{tikzcd}
 (O[k]\oplus  O[\ell]^{\circ})_0^{(+)} \arrow[r]\arrow[d] &  O[k+1+\ell]_0^{(+)}\arrow[d]\\
 O[0]_0^{(+)}\oplus  (O[0]^\circ)_0^{(+)}\arrow[r] & \Sigma( O[k]\Gray  O[\ell]^{\circ})_0^{(+)}.
\end{tikzcd}
    \end{equation}

Let $q>0$. Vandermonde's identity 
\[
{[k+1+\ell]\choose[q]}=
\sum_{i=-1}^{q}{[k]\choose[i]}\cdot{[\ell]\choose[q-1-i]}=
{[k]\choose[q]}+{[\ell]\choose[q]}+\sum_{i=0}^{q-1}{[k]\choose[i]}\cdot{[\ell]\choose[q-1-i]},
\]
can be equivalently expressed as a pushout of pointed sets
\[ \begin{tikzcd}
{[k]\choose [q]}_+\vee{[\ell]\choose[q]}_+ \arrow[r]\arrow[d] &{[k+1+\ell]\choose [q]}_+ \arrow[d]\\
 \{*\}\arrow[r] &  \bigvee\limits_{i=0}^{q-1}\left({[k]\choose[i]}\times{[\ell]\choose[q-1-i]}\right)_+,
\end{tikzcd}\]
Here, the cocomponents of the top horizontal map are the canonical, and the right vertical map is the one from \cref{descriptionphi}.
By \cref{ColimitsADCH}(2),
we then obtain a pushouts of abelian groups and of commutative monoids:
 \[
 \begin{tikzcd}
 \mathbb Z^{[k]\choose [q]}\oplus\mathbb Z^{[\ell]\choose[q]} \arrow[r]\arrow[d] &\mathbb Z^{[k+1+\ell]\choose [q]} \arrow[d]\\
 0\oplus0\arrow[r] &  \bigoplus\limits_{r=0}^q\mathbb Z^{[k]\choose[r]}\otimes\mathbb Z^{[\ell]\choose[q-1-r]}
\end{tikzcd}
\text{ and }
\begin{tikzcd}
 \mathbb N^{[k]\choose [q]}\oplus\mathbb N^{[\ell]\choose[q]} \arrow[r]\arrow[d] &\mathbb N^{[k+1+\ell]\choose [q]} \arrow[d]\\
 0\oplus0\arrow[r] &  \bigoplus\limits_{r=0}^q\mathbb N^{[k]\choose[r]}\otimes\mathbb N^{[\ell]\choose[q-1-r]}.
\end{tikzcd}
\]
This means that \eqref{sigmaoforientaltensor} induces a pushout of abelian groups (resp.~commutative monoids):
  \begin{equation} \label{sigmaoforientaltensorq}
\begin{tikzcd}
 (O[k]\oplus  O[\ell]^{\circ})_q^{(+)} \arrow[r]\arrow[d] &  O[k+1+\ell]_q^{(+)}\arrow[d]\\
 O[0]_q^{(+)}\oplus  (O[0]^\circ)_q^{(+)}\arrow[r] & \Sigma( O[k]\Gray  O[\ell]^{\circ})_q^{(+)}.
\end{tikzcd}
 \end{equation}
Combining \eqref{sigmaoforientaltensor0} and \eqref{sigmaoforientaltensorq}, by \cref{ColimitsADCH} we obtain
the desired pushout of augmented directed chain complexes \eqref{sigmaoforientaltensor}.
\end{proof}

\begin{prop}
\label{orientaltosuspensionepi}
Let $k,\ell\ge0$. The map
\[O[k+1+\ell]\to \chsusp\left(O[k]\otimes O[\ell]^{\circ}\right)\]
from \cref{orientaltosuspension} is an epimorphism of augmented directed chain complexes.
\end{prop}

\begin{proof}
Using the explicit computations from \cref{quotientoforiental} we see that for every $q\ge0$ the canonical map
\[ (O[k]\oplus  O[\ell]^{\circ})_q^{(+)} \to (O[0]\oplus  O[0]^{\circ})_q^{(+)}\]
is an epimorphism of abelian groups (resp.~commutative monoids). By \cref{ColimitsADCH2}
the canonical map
\[ O[k]\oplus  O[\ell]^{\circ} \to O[0]\oplus  O[0]^{\circ} \]
is then an epimorphism of augmented directed chain complexes. Given that epimorphisms are closed under pushout, by \cref{quotientoforiental}, the map
\[O[k+1+\ell]\to \chsusp\left(O[k]\otimes O[\ell]^{\circ}\right)\]
from \cref{orientaltosuspension} is then an epimorphism of augmented directed chain complexes, too.
\end{proof}

We can understand how to map orientals into suspensions:

\begin{prop}
For $m\ge1$ and $C$ an augmented directed chain complex, the map from \cref{orientaltosuspension} induces a natural bijection
\[
\coprod_{\substack{k+1+\ell=m,\\ k,\ell\geq-1}}\adch(O[k]\Gray O[\ell]^{\circ}, C) \stackrel{\cong}{\longrightarrow} \adch(O[m], \Sigma C).
\]
\end{prop}

\begin{proof}
For any $x\colon O[m]\to\Sigma C$ we set
\[k\coloneqq\#\{0\leq i\leq m\ |\ x([i])=\bot\}-1\text{ and }\ell\coloneqq m-1-k.\]
and construct a corresponding preimage $\widehat{x}\colon O[k]\Gray O[\ell]^{\circ}\to C$.

Let $k=-1$ (resp.~$k=m$). Then we take
\[\widehat x\colon O[-1]\Gray O[m] \cong O[-1]\to C,\text{ resp. }\widehat x\colon O[m]\Gray O[-1]\cong O[-1]\to C\] to be the trivial map.

Let $0\leq k\leq m-1$. By \cref{orientaltosuspensionepi}, the function
\[\Sigma\colon\coprod_{\substack{k+1+\ell=m\\0\leq k,\ell}}\adch(O[k]\Gray O[\ell]^{\circ}, C) \to
\coprod_{\substack{k+1+\ell=m\\ k,\ell\geq 0}}\prescript{{\Sigma O[-1]/}}{}{\adch}(\Sigma (O[k]\Gray O[\ell]^{\circ}), \Sigma C)\]
is bijective, so $x$ can be uniquely identified with the suspension map
\[\Sigma x\colon\Sigma(O[k]\Gray O[\ell]^{\circ})\to \Sigma C\]
under $\Sigma O[-1]$.
By composing with the map
\[\coprod_{\substack{k+1+\ell=m\\k,\ell\geq 0}}\prescript{{\Sigma O[-1]/}}{}{\adch}(\Sigma (O[k]\Gray O[\ell]^{\circ}), \Sigma C)\to
 \adch(O[m], \Sigma C)\]
induced by \cref{orientaltosuspension}, which is
is injective by \cref{suspensionFF}, $\Sigma x$ can be uniquely identified with a map of the form
\[\widehat x\colon O[k]\Gray O[\ell]^{\circ}\to C.\]
It is a straightforward verification that the assignment $x\mapsto\widehat x$ defines the inverse for the desired bijection.
\end{proof}

\section{\texorpdfstring{$\omega$}{w}-categories and algebraic models}

We recall the basic definitions around $\omega$-categories, as well as some constructions based on $\omega$-categories: the $\omega$-categorical suspension, the tensor product, the total dual, and Steiner's linearization, as well as the main properties that we use later in the paper, and relevant examples.

\subsection{\texorpdfstring{$\omega$}{w}-categories}

While we refer the reader to e.g.~\cite{StreetOrientedSimplexes} for a traditional approach to the definition of an $\omega$-category, we briefly recall the main features here.

The data of an \emph{$\omega$-category} $\cC$ consists of a collection of sets $\cC_q$, for ${q \geq 0}$,
where $\cC_0$ is called the set of \emph{objects} of $\cC$ and $\cC_q$ for $q>0$
is the set of \emph{$q$-cells} or cells of \emph{dimension} $q$ of $\cC$, together with:
\begin{itemize}[leftmargin=*]
    \item \emph{source} and \emph{target} operators $s_q, t_q \colon \cC_p \to \cC_q$
    for all $p > q \geq 0$;
    \item \emph{identity} operators $\id_q \colon \cC_p \to \cC_{q}$ for all $q\geq p\ge0$;
    \item \emph{composition} operators $\comp_p \colon \cC_q \times_{\cC_p} \cC_q \to \cC_q$
    defined for all $q > p \geq 0$ and all pairs of $q$-cells $(g, f)$ for which $s_p(g) = t_p f$.
\end{itemize}
We say that $\cC$ is an $\omega$-category if for all $r > q > p \geq 0$ the triple $(\cC_p, \cC_q, \cC_r)$
together with all the relevant source, target, identity and composition operators
is a $2$-category. In particular, 
\begin{equation}
\label{globularity}
s_ps_q f=s_p f    
\qquad\text{and}\qquad
t_pt_q f=t_p f    
\end{equation}
for any $r$-cell $f$ of $\cC$ and $r>q>p$.

An \emph{$\omega$-functor} $F \colon \cC \to \cD$ between $\omega$-categories $\cC$ and $\cD$ is a collection of maps $F_q \colon \cC_q \to \cD_q$ for
$q \geq 0$ that preserves source, target, identity, and composition operators.
We denote by $\omegacat$ the category of (small) $\omega$\nobreakdash-categories and $\omega$-functors.

A cell in an $\omega$-category $\cC$ is said to be \emph{trivial} if it is the identity of a cell of lower dimension. For $m\ge0$, an \emph{$m$-category} is an $\omega$-category in which all $q$-cells are trivial for $q>n$, and an \emph{$n$-functor} is an $\omega$-functor between $n$-categories. We denote by $n\cat$ the (full) subcategory of $\omega\cat$ given by $n$-categories and $n$-functors.

\begin{ex}
For $m\ge-1$, the \emph{$m$-oriental} $\cO[m]$ from \cite{StreetOrientedSimplexes}, \cite[Theorem 3.2]{SteinerOrientals} or \cite[\textsection7.2]{AraMaltsiniotisJoin}
is an $m$-category, and in particular an $\omega$-category. The reader who is not familiar with the original definition can also take the formula from \cref{orientalandnu} as the definition of $\cO[m]$.
\end{ex}


\begin{defn}[{\cite[\textsection1.8]{AraMaltsiniotisJoin}}]
\label{totaldual}
Let $\cC$ be an $\omega$-category.
The \emph{total dual}
of $\cC$
is the $\omega$-category $\cC^\circ$ with the following structure.
\begin{itemize}[leftmargin=*]
    \item The set of $q$-cells $\cC^\circ_q$ is
$\cC^\circ_q\coloneqq\cC_{q}$
\item The source map $s_{q}\colon\cC^\circ_{p}\to\cC^\circ_{q}$ is given by 
$s_q^{\cC^\circ} f=t_q^{\cC} f$  for all $p > q \geq 0$;
\item The target map $t_{q}\colon\cC^\circ_{p}\to\cC^\circ_{q}$ is given by 
$t_q^{\cC^\circ} f=s_q^{\cC} f$  for all $p > q \geq 0$;
\item The composition map $\ast_{p}\colon\cC^\circ_{q}\times_{\cC^\circ_{p}}\cC^\circ_{q}\to\cC^\circ_{q}$ is given by
$f\ast_{p}^{\cC^\circ} g=g \ast_{p}^{\cC} f$ for all $q > p \geq 0$
\item The identity map $\id_{q}\colon\cC^\circ_{p}\to\cC^\circ_{q}$ is given by $\id_q^{\cC^{\circ}} f=\id^{\cC}_q f$  for all $q\geq p\ge0$;
\end{itemize} 
The construction defines a functor $(-)^\circ\colon\omega\cat\to\omega\cat$.
\end{defn}

\subsection{Suspension of $\omega$-categories}

The following is a variant\footnote{Precisely, what we present in \cref{suspension} is the composite of the one used in \cite[\textsection B.6.5]{AraMaltsiniotisJoin} with the total dual from \cite[\textsection1.8]{AraMaltsiniotisJoin}.} of the construction treated in \cite[\textsection B.6.5]{AraMaltsiniotisJoin}. When the input is an ordinary $1$-category $\cC$, the suspension agrees with the one that we previously considered in \cite{ORfundamentalpushouts}.

\begin{defn}
\label{suspension}
Let $\cC$ be an $\omega$-category. 
The \emph{suspension}
of $\cC$
is the $\omega$-category $\Sigma\cC$ with the following structure.
\begin{itemize}[leftmargin=*]
    \item The set of $q$-cells $(\Sigma\cC)_q$ is
\[(\Sigma\cC)_q:=\{\bot,\top\}\cup\cC_{q-1},\quad\text{with}\quad(\Sigma\cC)_0:=\{\bot,\top\}\]
\item The source map $s_{q}\colon\Sigma\cC_p\to(\Sigma\cC)_{q}$ for $q>1$ is given by
\[s^{\Sigma\cC}_{q}f=s^{\cC}_{q-1} f,\quad s_{q}^{\Sigma\cC}\bot=\bot,\quad s_{q}^{\Sigma\cC}\top=\top,\quad\text{with}\quad s_0^{\Sigma\cC}f=\bot.\]
\item The target map $t_{q}\colon\Sigma\cC_p\to(\Sigma\cC)_{q}$ for $q>1$ is given by 
\[t^{\Sigma\cC}_{q}f=t^{\cC}_{q-1} f,\quad t_{q}^{\Sigma\cC}\bot=\bot,\quad t_{q}^{\Sigma\cC}\top=\top,\quad\text{with}\quad t_0^{\Sigma\cC}f=\top.\]
\item The identity map $\id_{q}\colon\Sigma\cC_{p}\to(\Sigma\cC)_{q}$ is given by \[\id_q^{\Sigma\cC} f=\id_q^{\cC} f,\quad \id_q^{\Sigma\cC}\bot=\bot, \quad \id_q^{\Sigma\cC}\top=\top.\]
\item The composition map $\ast_{p}\colon\Sigma\cC_{q}\times_{(\Sigma\cC)_p}(\Sigma\cC)_{q}\to(\Sigma\cC)_{q}$ is given by
\[g \ast_{p}^{\Sigma\cC} f = g \ast_{p-1}^{\cC} f\]
\end{itemize} 
Regarding $\Sigma \cC$ as an $\omega$-category bipointed on ${\bot}$ and ${\top}$,
the construction defines a functor $\Sigma\colon\omega\cat\to\omega\cat_{*,*}$.
\end{defn}


\subsection{Steiner's functors}

We briefly recall Steiner's adjoint pair that relates $\omega$-categories and augmented directed chain complexes. For a more detailed treatment, see \cite[Definition~2.8]{SteinerEmbedding} or~\cite[\textsection 2.4]{AraMaltsiniotisJoin}.

\begin{defn}
Let $C$ be an augmented directed chain complex. A \emph{Steiner table} in $C$ is a matrix
\[
	x=\tabld{x}{q}  
\]
such that, for $\alpha=+,-$ and $0\le p\le q$, the following hold:
\begin{enumerate}
	\item $x^\alpha_p$ belongs to $C^+_p$;
	\item $\partial(x^\alpha_p) = x^+_{p-1} - x^-_{p-1}$ for $0<p\leq q$;
	\item $\varepsilon(x^\alpha_0) = 1$;
	\item $x_q^- = x_q^+$.
\end{enumerate}
\end{defn}

\begin{defn}[{\cite[Definition~2.8]{SteinerEmbedding},~\cite[\textsection 2]{AraMaltsiniotisJoin}}]
Let $C$ be an augmented directed chain complex. The \emph{$\omega$-categorical realization} of $C$ is the $\omega$-category $\nu C$ is defined as follows.
\begin{itemize}[leftmargin=*]
    \item The set $(\nu C)_q$ of $q$-cells is given by
\[(\nu C)_q\coloneqq\left\{x=\tabld{x}{q}\ \right |\left.\ x\text{ Steiner table in }C \vphantom{\tabld{x}{q}}\right\}.\]
\item The source map $s_q\colon(\nu C)_{p}\to(\nu C)_q$ is given by
\[s_q\tabld{x}{p}\coloneqq\begin{pmatrix}x^-_0 &\dots &x^-_{q-1}
  &x^-_{q}\cr\noalign{\vskip 3pt} x^+_0 &\dots &x^+_{q-1}
  &x^-_{q}\end{pmatrix}\]
\item The target map $t_q\colon(\nu C)_{p}\to(\nu C)_q$ is given by
\[t_q\tabld{x}{p}\coloneqq\begin{pmatrix}x^-_0 &\dots &x^-_{q-1}
  &x^+_{q}\cr\noalign{\vskip 3pt}  x^+_0 &\dots &x^+_{q-1}
  &x^+_{q}\end{pmatrix}\]
\item The composition map $\ast_p\colon(\nu C)_{q}\times_{(\nu C)_p}(\nu C)_{q}\to(\nu C)_q$ is given by
\begin{align*}
&\tabld{x}{q}\ast_p\tabld{y}{q}\coloneqq\\
&\hphantom{XXXXXXX}\begin{pmatrix}x^-_0 &\dots &x^-_{p-1}
  &y^-_{p} & x_{p+1}^-+y_{p+1}^- &  \dots &x_{q}^-+y_{q}^-\cr\noalign{\vskip 3pt} 
  x^+_0 &\dots &x^+_{p-1}
  &x^+_{p} & x_{p+1}^++y_{p+1}^+ &  \dots &x_{q}^++y_{q}^+\end{pmatrix}
  \end{align*}
\item The identity map $\id_q\colon(\nu C)_{p}\to(\nu C)_q$ is given by
\[\id_q\tabld{x}{p}\coloneqq\underbrace{\begin{pmatrix}x^-_0 &\dots &x^-_{p-1}
  &x^-_{p} & 0 & 0 &0 \dots &0\cr\noalign{\vskip 3pt} 
  x^+_0 &\dots &x^+_{p-1}
  &x^+_{p} & 0 & 0 &0 \dots &0\end{pmatrix}}_{q+1}\]
\end{itemize}
The construction extends to a functor $\nu\colon\adch\to\omega\cat$.
\end{defn}

\begin{defn}
Let $\cC$ be an $\omega$-category.
The \emph{linearization} of $\cC$ is the augmented directed chain complex $\lambda\cC$ defined as follows.
\begin{itemize}[leftmargin=*]
    \item The abelian group $(\lambda\cC)_q$ of $q$-chains of $\lambda\cC$ is the quotient of $\mathbb Z[\cC_q]$ given by
\begin{equation}
\label{lambda}
(\lambda\cC)_q\coloneqq\frac{\mathbb Z[\cC_q]}{\left<[x\ast_p y]_q-[x]_q-[y]_q \ |\ x,y\in\cC_q;p<q\right>}.\ 
\end{equation}
\item The positivity submonoid $(\lambda\cC)_q^+$ is the submonoid of $(\lambda\cC)_q$ generated by the collection of elements $[f]_q$ for $f$ a $q$-cell of $\cC$.
\item The differential map $\partial_{q}\colon(\lambda\cC)_{q+1}\to(\lambda\cC)_{q}$ is determined by the condition on generators $f\in\cC_q$ given by
\[\partial_{q}([f]_{q+1})\coloneqq[t_{q} f]_{q}-[s_{q} f]_{q},\]
\item The augmentation map $\varepsilon\colon(\lambda\cC)_0\to\mathbb Z$ is determined by the condition on generators $x\in\cC_0$ given by
\[\varepsilon([x]_0)\coloneqq1.\]
\end{itemize}
The construction extends to a functor $\lambda\colon\omega\cat\to\adch$.
\end{defn}

\begin{thm}[{\cite[\textsection2]{SteinerEmbedding}}]
\label{nulamba}
The functors $\nu$ and $\lambda$ form an adjoint pair
\[\lambda\colon \omega\cat\rightleftarrows\adch\colon\nu.\]
In other words, for any $\omega$-category $\cC$ and any augmented directed chain complex $\overline C$ there is a natural bijection
\[\adch(\lambda \cC,\overline C)\cong\omega\cat(\cC,\nu \overline C).\]
\end{thm}

\begin{ex}[{\cite[Theorem~3.2]{SteinerOrientals}}]
\label{orientalandnu}
For $m\ge0$, there is an isomorphism of augmented directed chain complexes
\[\lambda \cO[m]\cong O[m]\]
and an isomorphism of $\omega$-categories
\[\cO[m]\cong\nu O[m].\]
\end{ex}

\begin{lem}[{\cite[Proposition 2.19]{AraMaltsiniotisJoin}}]
\label{opandnu}
Let $C$ be an augmented directed chain complex and $\cC$ an $\omega$-category.
\begin{enumerate}[leftmargin=*]
    \item There is a natural isomorphism of $\omega$-categories
\[\nu( C^\circ) \cong  (\nu C)^\circ.\]
    \item There is a natural isomorphism of augmented directed chain complexes
\[\lambda( \cC^\circ) \cong  (\lambda \cC)^\circ.\]
\end{enumerate}
\end{lem}

\subsection{Steiner's functors and suspension}

\begin{lem}
\label{suspensionandnu}
For an augmented chain complex $C$, there is a natural isomorphism
\[\nu\Sigma C \cong \Sigma \nu C.\]
\end{lem}

\begin{proof}
The prototypical element of both $(\nu\Sigma C)_q$ and $(\Sigma\nu C)_q$ is can be expressed as a table of the form
\[x=\begin{pmatrix}\bot&x^-_0 &\dots &x^-_{q-2}
  &x^-_{q-1}\cr\noalign{\vskip 3pt} \top & x^+_0 &\dots &x^+_{q-2}
  &x^+_{q-1}\end{pmatrix},\]
  where
  \begin{enumerate}
	\item $x^\alpha_p$ belongs to $C^+_p$;
	\item $\partial(x^\alpha_p) = x^+_{p-1} - x^-_{p-1}$ for $0<p< q$;
	\item $\varepsilon(x^\alpha_0) = 1$;
	\item $x_{q-1}^- = x_{q-1}^+$.
\end{enumerate}
One can check that this identification is compatible with source, target, identity and composition operations, and the desired isomorphism of $\omega$-categories follows.
\end{proof}

\subsection{Tensor product of $\omega$-categories}

The statement of the following theorem relies on the notion of a \emph{strong Steiner $\omega$-category}, a.k.a. \emph{$\omega$-category that admits a strongly loop-free atomic basis}. Those are particularly nice $\omega$-categories that are in a sense ``free'' and ``loop-free'' and we refer the reader to \cite{SteinerEmbedding,AraMaltsiniotisJoin,AGOR} for an account on strong Steiner $\omega$-categories. There is also a notion of a \emph{strong Steiner complex}, a.k.a.~\emph{augmented directed chain complex that admits a strongly loop-free and unital basis}, which correspond in a precise sense to Strong Steiner $\omega$-categories under the adjunction $(\lambda,\nu)$. For the purpose of this paper, it is sufficient to know the following.
\begin{itemize}[leftmargin=*]
    \item For every $m\ge0$ the $m$-oriental $\cO[m]$ is a strong Steiner $\omega$-category (as shown in \cite[Example 3.8]{SteinerEmbedding}), and so it its total dual $\cO[m]^\circ$ (which can be verified directly).
    \item For every $m\ge0$ the $m$-cell is a strong Steiner $\omega$-category (as shown in \cite[Example 3.9]{SteinerEmbedding}).
    \item
   For any strong Steiner $\omega$-category  $\cC$, the unit of the adjunction from \cref{nulamba} is an isomorphism of $\omega$-categories $\eta_\cC\colon\cC\cong\nu\lambda\cC$ (as shown in \cite[Theorem 5.11]{SteinerEmbedding}).
   \item For any strong Steiner complex $C$, the counit of the adjunction from \cref{nulamba} is an isomorphism of augmented directed chain complexes $\epsilon_C\colon\lambda\nu C\cong C$ (as shown in \cite[Theorem 5.11]{SteinerEmbedding}).
   \item For any strong Steiner $\omega$-category $\cC$, the augmented directed chain complex $\lambda\cC$ is a strong Steiner complex (as shown in \cite[Theorem 5.11]{SteinerEmbedding}).
\item For any strong Steiner complex $C$, the $\omega$-category $\nu C$ is a strong Steiner $\omega$-category (as shown in \cite[Theorem 5.11]{SteinerEmbedding}).
    \item For any strong Steiner complex $C$ and $\overline{\cC}$ the augmented directed chain complex $C\otimes\overline C$ is a strong Steiner complex (as shown in \cite[Proposition A.3]{AraMaltsiniotisJoin}).
\item for any strong Steiner complexes $C$ and $\overline C$
there is a natural isomorphism of augmented directed chain complexes
$\nu C\otimes\nu \overline C\cong\nu(C\otimes \overline C)$ (as shown in \cite[Theorem A.15]{AraMaltsiniotisJoin}).
\end{itemize}

\begin{thm}[{\cite[Theorem A.15]{AraMaltsiniotisJoin}}]
\label{tensorcompatible}
There exists a unique -- up to unique monoidal isomorphism -- monoidal structure $\otimes\colon\omega\cat\times\omega\cat\to\omega\cat$ on $\omegacat$, called the \emph{tensor product} of $\omega$-categories, such that
\begin{itemize}[leftmargin=*]
\item for any strong Steiner $\omega$-categories $\cC$ and $\overline{\cC}$ the tensor product  $\cC\otimes\overline{\cC}$ is the $\omega$-category
\[\cC\otimes\overline{\cC}\coloneqq\nu(\lambda\cC\otimes\lambda\overline{\cC});\]
\item the functor $-\otimes-$ commutes with colimits in each variable.
\end{itemize}
\end{thm}

\begin{prop}
\label{lambdamonoidal}
The linearization functor defines a strong monoidal functor $\lambda\colon(\omega\cat,\otimes)\to(\adch,\otimes)$. That is, for any $\omega$-categories $\cC$ and $\overline\cC$ there is a natural isomorphism of augmented directed chain complexes
\[\lambda(\cC\otimes\overline\cC)\cong\lambda\cC\otimes\lambda\overline\cC.\]
\end{prop}

\begin{proof}
First, we observe that for any strong Steiner $\omega$-categories $\cC$ and $\overline{\cC}$ we have
\[
\begin{array}{llll}
     \lambda(\cC\otimes\overline\cC)& \cong \lambda(\nu\lambda\cC\otimes\nu\lambda\overline\cC)
     \\
     & \cong \lambda\nu(\lambda\cC\otimes\lambda\overline\cC)
     \\
     & \cong \lambda\cC\otimes\lambda\overline\cC. 
\end{array}
\]
so the desired isomorphism holds for strong Steiner $\omega$-categories.
Since any $\omega$-category is a colimit of strong Steiner $\omega$-categories (as cells are in particular strong Steiner $\omega$-categories), and the functors $\lambda(-\otimes-),(\lambda-)\otimes(\lambda-)\colon\omega\cat\times\omega\cat\to\omega\cat$ commute with colimits in both variables, the desired isomorphism follows.
\end{proof}

We can understand tensor product of orientals:

\begin{ex}
\label{tensoroforientals}
Let $k,\ell\geq -1$. 
Have isomorphism of $\omega$-categories
\[
\begin{array}{llll}
    \cO[k]\otimes\cO[\ell]^\circ   & \cong& \nu O[k]\Gray (\nu O[\ell])^{\circ}&\text{\cref{orientalandnu}}\\
    & \cong& \nu O[k]\Gray \nu( O[\ell]^{\circ})&\text{\cref{opandnu}(1)}\\
    & \cong& \nu (O[k]\Gray O[\ell]^{\circ})&\text{\cref{tensorcompatible}}\\ 
\end{array}\]
and of augmented directed chain complexes
\[
\begin{array}{llll}
   \lambda(\cO[k]\otimes\cO[\ell]^\circ)  &\cong  & \lambda\cO[k]\otimes\lambda(\cO[\ell]^{\circ})&\text{\cref{lambdamonoidal}}\\
   &\cong  & \lambda\cO[k]\otimes(\lambda\cO[\ell])^{\circ}&\text{\cref{opandnu}(2)}\\
    &\cong&O[k]\Gray O[\ell]^{\circ}&\text{\cref{orientalandnu}}
\end{array}\]
\end{ex}

Can understand the suspension of tensor product of orientals:

\begin{rmk}
\label{orientaltosuspension2}
Let $k,\ell\geq0$.
Applying $\nu$ to the square \eqref{sigmaoforientaltensor} -- and evoking \cref{tensoroforientals,orientalandnu,suspensionandnu}) -- we obtain the diagram of $\omega$-categories
    \begin{equation}
\label{sigmaoforientaltensor2}
\begin{tikzcd}
 \cO[k]\oplus  \cO[\ell]^{\circ} \arrow[r]\arrow[d] &  \cO[k+1+\ell]\arrow[d]\\
 \cO[0]\oplus  \cO[0]\arrow[r] & \Sigma( \cO[k]\Gray  \cO[\ell]^{\circ}).
\end{tikzcd}
    \end{equation}
 In particular, the map
\begin{equation}
    \label{orientaltosuspension2}
\cO[k+1+\ell]\to\Sigma( \cO[k]\Gray  \cO[\ell]^{\circ}).
\end{equation}
is induced by \eqref{orientaltosuspension}.
\end{rmk}

\begin{prop}
 \label{quotientoforiental2} Let $k,\ell\geq -1$. The diagram \eqref{sigmaoforientaltensor2} induces a natural isomorphism of $\omega$-categories
 \[\Sigma( \cO[k]\Gray  \cO[\ell]^{\circ})\cong(\cO[0]\oplus  \cO[0]^{\circ})\aamalg{\cO[k]\oplus  \cO[\ell]^\circ}\cO[k+1+\ell].\]
\end{prop}

\begin{proof}
Consider the commutative diagram of augmented directed chain complexes on the left, and the induced commutative diagram of $\omega$-categories on the right:
\[
\begin{tikzcd}
O[k]\oplus  O[\ell]^\circ\arrow[r]\arrow[d]&O[k+1+\ell]\arrow[d]\\O[0]\oplus  O[0]^{\circ}\arrow[r]&\Sigma( O[k]\Gray  O[\ell]^{\circ})
\end{tikzcd}
\quad
\quad
\rightsquigarrow
\quad
\quad
\begin{tikzcd}
\cO[k]\oplus  \cO[\ell]^\circ\arrow[r]\arrow[d]&\cO[k+1+\ell]\arrow[d]\\\cO[0]\oplus  \cO[0]^{\circ}\arrow[r]&\Sigma( \cO[k]\Gray  \cO[\ell]^{\circ}).
\end{tikzcd}
\]
The square on the left is a pushout by \eqref{sigmaoforientaltensor} and, as an application of \cite[Th\'eor\`eme~3.1.5]{LoubatonNerfs},
so is the pushout on the right.
\end{proof}

Can understand how to map orientals into suspension $\omega$-categories:

\begin{prop}
\label{orientaltosuspension2}
For $m\ge1$ and $\cC$ an $\omega$-category of the form $\cC\cong\nu C$, there is a natural bijection
\[
\coprod_{\substack {k+1+\ell=m\\ k,\ell\geq-1}}\omega\cat(\cO[k]\Gray \cO[\ell]^{\circ}, \cC) \stackrel{\cong}{\longrightarrow}\omega\cat(\cO[m], \Sigma \cC).
\]
\end{prop}

\begin{proof}
There's a natural bijection
\[
\begin{array}{lllll}
  \coprod\limits_{k,\ell}\omega\cat(\cO[k]\Gray \cO[\ell]^{\circ}, \cC)  & \cong &    \coprod\limits_{k,\ell}\omega\cat(\cO[k]\Gray \cO[\ell]^{\circ}, \nu C)\\
& \cong &    \coprod\limits_{k,\ell}\adch(\lambda(\cO[k]\Gray \cO[\ell]^{\circ}) ,C)&\text{\cref{nulamba}}\\
& \cong &    \coprod\limits_{k,\ell}\adch(O[k]\Gray O[\ell]^{\circ},C)&\text{\cref{orientalandnu}}\\
\end{array}
\]
and a natural bijection
\[\begin{array}{lllll}
\omega\cat(\cO[m], \Sigma\cC)&\cong&\omega\cat(\cO[m], \Sigma\nu C)&\text{}\\
&\cong&\omega\cat(\cO[m], \nu\Sigma C)&\text{\cref{suspensionandnu}}\\
&\cong&\adch(\lambda\cO[m], \Sigma C)&\text{\cref{nulamba}}\\
&\cong&\adch(O[m], \Sigma C)&\text{\cref{orientalandnu}}\\
\end{array}
\]

They fit into a commutative diagram of sets
\[\begin{tikzcd}[column sep=0.25cm]
\omega\cat(\cO[m],\Sigma\cC)\arrow[r,<-,""]\arrow[d,<- ,"\cong"]&
\coprod\limits_{k,\ell}\omega\cat(\Sigma(\cO[k]\Gray \cO[\ell]^{\circ}),\Sigma\cC)\arrow[d,<- ,"" swap]&\coprod\limits_{k,\ell}\omega\cat(\cO[k]\Gray \cO[\ell]^{\circ}, \cC)\arrow[l,"\Sigma"]\arrow[d,<- ,"\cong" swap] \arrow[ll, bend right=15]\\
\adch(O[m],\Sigma C)\arrow[r,<-,""]&
\coprod\limits_{k,\ell}\adch(\Sigma(O[k]\Gray O[\ell]^{\circ}),\Sigma C)&\coprod\limits_{k,\ell}\adch(O[k]\Gray O[\ell]^{\circ}, C)\arrow[l,"\Sigma"]\arrow[ll, bend left=20]
\end{tikzcd}\]
so we obtain two equal bijections of the desired form.
\end{proof}


\section{Complicial sets and complicial nerve of $\omega$-categories}

We recall the basic definitions around simplicial sets with marking and $n$-complicial sets, as well as some constructions based on simplicial sets with marking: the suspension and the complicial nerve, as well as the main properties that we use later in the paper, and relevant examples. The study of the homotopy theory of complicial sets originated with \cite{VerityComplicialI,VeritySlides}, and continued with \cite{RiehlNotes,ORfundamentalpushouts,OR,RiehlVerityBook}.

\subsection{Complicial sets}

We recall the main facts about complicial sets that will be used in this paper. 
 
\begin{defn}
A \emph{simplicial set with marking} is a pair $(X,tX)$ where $X$ is a simplicial set, and $tX=\coprod_{m\ge1}tX_m\subseteq\coprod_{m\ge1}X_m$ is a collection of subsets of simplices of $X$, called \emph{marked} simplices, which have positive dimension that contain all degenerate simplices of $X$.
\end{defn}

We denote by $m\sset$ the category of simplicial sets with marking and marking-preserving simplicial maps.

\begin{rmk}
The category $\msset$ is cocomplete, and colimits are computed degreewise (a simplex is marked in a colimit if it admits a marked representative).
\end{rmk}

\begin{defn}
A sub-simplicial set with marking $X$ of a simplicial set with marking $Y$ is \emph{regular} if a simplex of $X$ is marked in $X$ precisely when it is marked in $Y$.
\end{defn}

\begin{notn}
\label{preliminarynotation}
We denote
\begin{itemize}[leftmargin=*]
    \item by $\Delta^k[m]$, for $0\leq k \leq m$, the standard $m$-simplex in which a non-degenerate simplex is marked if and only if it contains the vertices $\{k-1,k,k+1\}\cap [m]$;
    \item by $\Delta^k[m]'$, for $0\leq k \leq m$, the standard $m$-simplex with marking obtained from $\Delta^k[m]$ by additionally marking the $(k-1)$-st and $(k+1)$-st face of $\Delta[m]$;
    \item by $\Delta^k[m]''$, for $0\leq k \leq m$, the standard $m$-simplex with marking obtained from $\Delta^k[m]'$ by additionally marking the $k$-th face of $\Delta[m]$;
    \item by $\Lambda^k[m]$, for $0\leq k \leq m$, the regular sub-simplicial set of $\Delta^k[m]$ with marking whose simplicial set is the $k$-horn $\Lambda^k[m]$.
    \item by $\Delta[3]_{\sharp}$ the standard $3$-simplex with the maximal marking.
    \item by $\Delta[3]_{\mathrm eq}$ the standard $3$-simplex in which the $1$-simplices $[0,2]$ and $[1,3]$ are marked, as well as all simplices in dimension $2$ or higher.
\end{itemize}
\end{notn}

The following class of maps plays a role in the model structure on $\msset$ for $(\infty,n)$-categories, with $n\in\mathbb N\cup\{\infty\}$.

\begin{defn}
\label{complicialextensions}
Let $n\in\mathbb N\cup\{\infty\}$.
\begin{enumerate}[leftmargin=*]
\item For $m> 1$ and $0< k< m$, the \emph{complicial inner horn extension} is the inclusion
\[\Lambda^k[m]\to \Delta^k[m].\]
 \item For $m\geq 2$ and $0< k < m$, the \emph{complicial thinness extension} is the inclusion
\[\Delta^k[m]' \to \Delta^k[m]''.\]
\item For $m>n$, the \emph{triviality extension} is the inclusion
\[\Delta[m]\to\Delta[m]_t.\]
\item For $m\ge-1$, the \emph{complicial saturation extension} is the inclusion\footnote{The reader can find the join of simplicial sets with marking $\star\colon\msset\times\msset\to\msset$ in \cite[\textsection3.1]{VerityComplicialI}, but it will not be needed explicitly in this paper.}
\[\Delta[3]_{\mathrm{eq}}\star\Delta[m]\to\Delta[3]_{\sharp}\star\Delta[m].\]
\end{enumerate}
\end{defn}

We fix the following terminology (cf.~\cite[Def.~15]{VerityComplicialI}).

\begin{defn}
A map of simplicial sets with marking $X\to Y$ is a \emph{complicial inner anodyne extension} if it can be written as a retract of a transfinite composition of pushouts of maps of type (1) and (2) from \cref{complicialextensions}.
\end{defn}

\begin{rmk}
 \label{underlyingcomplicialinner}
One can prove with standard model categorical techniques the following formal properties of complicial inner anodyne extensions.
\begin{enumerate}[leftmargin=*]
\item The underlying simplicial map of a complicial inner anodyne extension is an inner anodyne extension of simplicial sets.
     \item The class of complicial inner anodyne extensions is closed under transfinite composition and pushouts.
\end{enumerate}
\end{rmk}

\begin{lem}[{\cite[Lemma~1.12]{ORfundamentalpushouts}}]
  \label{CompMarkAtOnce} For $m\geq 2$ and $0<k<m$, let $\Lambda^k[m]'$ denote the regular subset of $\Delta^k[m]'$ whose underlying simplicial set is given by the $k$-horn $\Lambda^k[m]$. The inclusion
    \[\Lambda^k[m]'\to \Delta^k[m]''\]
is a complicial inner anodyne extension.
\end{lem}

The category $\msset$ hosts a model for $(\infty,n)$-categories.

\begin{defn}
\label{defcomplicial}
An \emph{$n$-complicial set}
is a simplicial set that has the right lifting property with respect to all maps of type (1)-(4) from \cref{complicialextensions}.
\end{defn}

\begin{thm}[{\cite[Theorem~1.25]{OR}}]
\label{modelstructurewithsaturation}
Let $n\in\mathbb N\cup\{\infty\}$.
\label{modelstructureondiscretepresheaves}
The category $m\sset$ supports a
cartesian closed model structure $\msset_{(\infty,n)}$, that we call the \emph{model structure for $(\infty,n)$-categories}, where \begin{itemize}[leftmargin=*]
    \item the fibrant objects are precisely the \emph{saturated $n$-complicial sets},
    \item the cofibrations are precisely the monomorphisms \textnormal{(}of underlying simplicial sets\textnormal{)},
    \item all complicial anodyne extensions are weak equivalences.
\end{itemize} 
\end{thm}

\begin{rmk}
\label{nonsaturatedmodelstructures}
Other model structures are sometimes considered on $\msset$, for instance those from
\cite[Theorem 100]{VerityComplicialI} and\cite[Examples~3.33--3.36]{RiehlNotes}. But in all the aforementioned model structures complicial inner anodyne extensions are weak equivalences.
\end{rmk}


\subsection{Suspension of complicial sets}

We now define the suspension
of simplicial sets with marking.

 \begin{defn}[{\cite[Definition~2.6]{ORfundamentalpushouts}}]
 \label{DefSuspension}
 Let $X$ be a simplicial set with marking.
The \emph{suspension} $\Sigma X$ of $X$ is the simplicial set with marking defined as follows.
\begin{itemize}[leftmargin=*]
    \item The set $(\Sigma X)_m$ of $m$-simplices is given by
\[(\Sigma X)_{m}=\left\{
\begin{array}{llll}
  \{\bot,\top\}   & \text{if $m=0$,} \\
  \{\bot,\top\}\cup\coprod\limits_{k=0}^{m-1}X_{k}   & \text{if $m=0$.}
\end{array}
\right.\]
\item The face map $d_i\colon(\Sigma X)_m\to(\Sigma X)_{m-1}$ satisfies
\[d_i\bot=\bot\text{ and }d_i\top=\top,\]
and restricts to the map $X_k\subseteq(\Sigma X)_m\to(\Sigma X)_{m-1}$ given by
\[
d_i(x)=\left\{\begin{array}{lllll}
d_ix\in X_{k-1}\subseteq (\Sigma X)_{m-1}&\mbox{ if }0\leq i\leq k,\\
x\in X_{k}\subseteq (\Sigma X)_{m} &\mbox{ if }k+1\leq i \leq m,
\end{array}
\right.
\]
\item The degeneracy map $s_i\colon(\Sigma X)_{m}\to(\Sigma X)_{m+1}$ 
satisfies
\[s_i\bot=\bot\text{ and }s_i\top=\top,\]
and
restricts to the map $X_k\subseteq(\Sigma X)_m\to(\Sigma X)_{m+1}$ given by
\[
s_i(x)=
\left\{\begin{array}{ll}
s_ix\in X_{k+1}\subseteq (\Sigma X)_{m+1} & \mbox{ if }0\leq i\leq k,\\
x\in X_{k}\subseteq (\Sigma X)_{m}& \mbox{ if }k+1\leq i \leq m.
\end{array}
\right.
\]
\item The set $t(\Sigma X)_m$ of marked $k$-simplices is given by
\[t(\Sigma X)_{m}=\coprod_{k=0}^{m}tX_{k}.\] 
\end{itemize}
Regarding $\Sigma X$ as a simplicial set with marking bipointed on ${\bot}$ and ${\top}$, the construction defines a functor $\Sigma\colon\msset\to\msset_{*,*}$.
 \end{defn}

\begin{rmk}
The set $(\Sigma X)_m^{\textrm{nd}}$ of non-degenerate $m$-simplices of $\Sigma X$ for $m>0$ is contained in the set of the non-degenerate $(m-1)$-simplices of $X$, namely
\[(\Sigma X)_m^{\textrm{nd}}\subseteq X_{m-1}.\]
\end{rmk}

As a special case of the slice model structures, constructed e.g.~in \cite{HirschhornOvercategories}, we also obtain model structure $(\msset_{(\infty,n+1)})_{*,*}$ on the category $m\sset_{*,*}$ of bi-pointed simplicial sets with marking.

\begin{lem}[{\cite[Lemma 2.7]{ORfundamentalpushouts}}]
\label{suspensionhomotopical}
The marked suspension defines a left Quillen functor
\[\Sigma\colon\msset_{(\infty,n)}\to (\msset_{(\infty,n+1)})_{*,*}.\]
\end{lem}


\subsection{Complicial nerve of \texorpdfstring{$\omega$}{w}-categories}

The geometry of orientals is such that the construction $m\mapsto \cO[m]$ defines a cosimplicial object $\cO[\bullet]$ in $\omega\cat$, and in particular it makes sense to define the nerve construction $N\colon\omega\cat\to\sset$ originally due to Street {\cite{StreetOrientedSimplexes}}.
The Street nerve can be endowed with the following marking, originally considered by Roberts in unpublished work and Street in \cite{StreetOrientedSimplexes}, further studied by Verity in \cite{VerityComplicialAMS}, and later discussed by Riehl in \cite{RiehlNotes}, obtaining a functor $\NRS\colon\omega\cat\to\msset$.

\begin{defn}
\label{RSnerve}
Let $\cC$ be an $\omega$-category. The \emph{Roberts--Street nerve} of $\cC$ is the simplicial set with marking defined as follows:
\begin{itemize}[leftmargin=*]
    \item The set of $m$-simplices is the set of $\omega$-functors $\mathcal O[m]\to \cC$, namely
\[N_m\cC=\omega\cat(\cO[m],\cC),\]
\item  the simplicial structure is induced by the geometry of orientals.
\item an $m$-simplex of $N\cC$ is marked in $\NRS\cC$ if and only if the corresponding $\omega$-functor $\cO[m]\to\cC$ sends the unique non-trivial $m$-cell $\langle[0,1,\dots,m]\rangle$ of $\cO[m]$ to a trivial $m$-cell of $\cC$, namely
\[x\in t(N\cC)_m\quad\Leftrightarrow\quad x(\langle[0,1,\dots,m]\rangle)=\id,\]
where $\langle[0,1,\dots,m]\rangle$ denotes the top non-identity $m$-cell of $\cO[m]$.
\end{itemize}
The construction extends to a functor $\NRS\colon\omega\cat\to\msset$.
\end{defn}
 In particular, in the Street nerve of an $n$-category $\cC$ all simplices in dimension at least $n+1$ are marked.

    \subsection{Complicial nerve and suspension}

We describe a comparison map between $\Sigma\NRS\cC$ and $\NRS\Sigma\cC$, which we will show to be furthermore a weak equivalence for an $\omega$-category of the form $\cC \cong \nu C$.

The simplicial sets have the same sets of $0$-simplices, namely
\[(\Sigma\NRS\cC)_0=\{\bot,\top\}=(\NRS\Sigma\cC)_0,\]
and we now analyze the set of $m$-simplices for $m>0$.

\begin{rmk}
\label{nerveofsuspension}
We have the following description for the set $(N\Sigma\cC)_m$ of $m$-simplices of $N\Sigma\cC$:
\[(N\Sigma\cC)_m=\omega\cat(\cO[m],\Sigma\cC).\]
Moreover,
\[x\in t(N\Sigma\cC)_m\quad\quad\Leftrightarrow\quad\quad x(\langle[0,1,\dots,m]\rangle)=\id.\]
If furthermore $\cC\cong\nu C$, by \cref{nulamba,orientaltosuspension,orientaltosuspension2} we also obtain the ``algebraic'' descriptions of $(N\Sigma\cC)_m$:
\[\begin{tikzcd}
(N\Sigma\cC)_m\arrow[r,equal]&[-0.6cm]\omega\cat(\cO[m],\Sigma\cC)\arrow[r,"\cong"]\arrow[d,<- ,"\cong"]&
\coprod\limits_{\substack{k+1+\ell=m,\\ k, \ell\geq -1} }\omega\cat(\cO[k]\Gray \cO[\ell]^{\circ}, \cC)\arrow[d,<- ,"\cong" swap]\\
&\adch(O[m],\Sigma C)\arrow[r,"\cong"]&
\coprod\limits_{\substack{k+1+\ell=m,\\ k, \ell\geq -1}}\adch(O[k]\Gray O[\ell]^{\circ}, C)
\end{tikzcd}\]
Moreover, 
\[x\in t(N\Sigma\cC)_m\quad\quad\Leftrightarrow\quad\quad x([ 0,\dots,k]\otimes[0,\dots \ell])=0.\]
\end{rmk}

Recall the bijection from \cref{orientaltosuspension2} which we use in the following definition.
\begin{defn}
\label{definitiontype}
Let $\cC$ be an $\omega$-category and $x\in(N\Sigma\cC)_m$. We say that the $m$-simplex $x$ if of \emph{type} $k$ if
\[x\in\omega\cat(\cO[k]\Gray \cO[\ell]^{\circ},\cC)\subseteq(N\Sigma\cC)_m.\]
\end{defn}

\begin{rmk}
\label{descriptionsimplices}
Let $\cC$ be an $\omega$-category and $x\in(N\Sigma\cC)_m$.
The following are equivalent.
\begin{enumerate}[leftmargin=*]
\item The simplex $x$ is the degeneracy of a $0$-simplex, namely $x=s_0^m\bot$ or $x=s_0^m\bot$.
\item The simplex $x$ has type $-1$ or $m$.
\end{enumerate}
If the equivalent conditions are met, we say that $x$ is \emph{totally degenerate}.
\end{rmk}

\begin{rmk}
\label{suspensionofnerve}
    Combining results from previous sections, we have the following equivalent descriptions for the set of $m$-simplices of $\Sigma N\cC$:
    \[ (\Sigma N\cC)_m\cong \{s_0^m\bot,s_0^m\top\}\amalg \coprod\limits_{k=0}^{m-1}(N\cC)_{k}\cong\{s_0^m\bot,s_0^m\top\}\amalg \coprod\limits_{k=0}^{m-1}\omega\cat(\cO[k],\cC).\]
    Moreover, for a non-totally degenerate simplex $x$, we have
\[x\in t(\Sigma N_m\cC)\quad\quad\Leftrightarrow\quad\quad  k<m-1\mbox{ or }x(\langle[0,1,\dots,k]\rangle)=\id.\] 
If furthermore $\cC\cong\nu C$, by \cref{nulamba} also get the ``algebraic'' descriptions:
  \[
  \begin{array}{llll}
       (\Sigma N\cC)_m\cong \{s_0^m\bot,s_0^m\top\}\amalg \coprod\limits_{k=0}^{m-1}(N\cC)_{k}&\cong\{s_0^m\bot,s_0^m\top\}\amalg \coprod\limits_{k=0}^{m-1}\omega\cat(\cO[k],\cC)  \\
       &\cong\{s_0^m\bot,s_0^m\top\}\amalg \coprod\limits_{k=0}^{m-1}\adch(O[k],C).
  \end{array}
  \]
Moreover,
\[x\in t(\Sigma N\cC)_m\quad\quad\Leftrightarrow\quad\quad k<m-1\mbox{ or }x([0,\dots,k])=0.\]
  \end{rmk}

  The canonical map(s) from either \eqref{orientaltosuspension1} or \eqref{sigmaoforientaltensor2} then induce a canonical natural map $(\Sigma N\cC)_m\to (N\Sigma\cC)_m$ which assembles into a map
 $\Sigma N\cC\to N\Sigma\cC.$

  \begin{prop}
 \label{comparisonmap}
For any $\omega$-category of the form $\cC\cong\nu C$, either of the maps \cref{orientaltosuspension} or \cref{orientaltosuspension2} induces:
\begin{enumerate}[leftmargin=*, label=(\arabic*)]
    \item\label{NonmarkedIncl} a natural inclusion of simplicial sets
 \[\Sigma N\cC\to N\Sigma\cC;\]
    \item\label{MarkedIncl} a natural regular inclusion of simplicial sets with marking
 \[\Sigma\NRS\cC\to \NRS\Sigma\cC.\]
\end{enumerate}
\end{prop}
 
 \begin{proof}
 At the level of simplicial sets, both inclusions act as follows:
  \begin{itemize}[leftmargin=*]
    \item they are identities on $0$-simplices, namely they send the $0$-simplex ${\bot}$ to ${\bot}$ and the $0$-simplex ${\top}$ to ${\top}$, and
\item referring to the identifications \cref{nerveofsuspension,suspensionofnerve}, they act on
an $(m+1)$-simplex $y\colon\cO[k]\to\cC$ of $\Sigma (N\cC)$ as
\[
[y\colon\cO[k]\to\cC] \mapsto [\cO[m] \to \Sigma(\cO[k] \otimes \cO[\ell]^{\circ}) \to \Sigma(\cO[k] \otimes \cO[0])\cong\Sigma(\cO[k])  \xrightarrow{\Sigma y } \Sigma \cC].
\]
 \end{itemize}
From this explicit description, we see that both maps are inclusions, and that the second one is regular, namely, a non-degenerate $(m+1)$-simplex of $\Sigma N\cC$ is marked in $\Sigma\NRS\cC$ if and only if the corresponding $m$-simplex of $N\cC$ is marked in $\NRS\cC$.
  \end{proof}

 \begin{rmk}
Given an $\omega$-category of the form $\cC\cong\nu C$, the map from \cref{comparisonmap} can be seen as induced by the canonical map
\[\NRS\cC\star\Delta[0]\to \NRS(\cC\star\Delta[0])\]
from \cite[Theorem~5.2]{GOR1}.
\[
\begin{tikzcd}
\NRS\cC\star\Delta[0]\arrow[r]\arrow[d]&\NRS(\cC\star\Delta[0])\arrow[d]\\
\Sigma\NRS\cC\arrow[r,dashed]&\NRS\Sigma\cC
\end{tikzcd}
\]
 \end{rmk}
 
We now prove that the comparison map is a weak equivalence if $\cC\cong\nu C$.

 \begin{thm}
 \label{ThmB}
Let $\cC$ be an $\omega$-category of the form $\cC\cong\nu C$.
\begin{enumerate}[leftmargin=*]
    \item \label{ThmBpart1}
    The inclusion from \cref{comparisonmap}\ref{NonmarkedIncl} is an inner anodyne extension, and in particular a weak equivalence in the Joyal model structure $\sset_{(\infty,1)}$
  \[\Sigma N\cC\xhookrightarrow{\hphantom{a}\simeq\hphantom{a}} N\Sigma\cC.\]
\item \label{ThmBpart2} The inclusion from \cref{comparisonmap}\ref{NonmarkedIncl} is a complicial inner anodyne extension, and in particular a weak equivalence in the model structure\footnote{The same map is also a weak equivalences in the model structures mentioned in \cref{nonsaturatedmodelstructures}.} $\msset_{(\infty,n)}$.
  \[\Sigma \NRS\cC\xhookrightarrow{\hphantom{a}\simeq\hphantom{a}} \NRS\Sigma\cC.\]
\end{enumerate}
 \end{thm}

The proof is given in the coming subsections.

\subsection{Proof of \cref{ThmB}}
Given an $\omega$-category $\cC\cong\nu C$, by \cref{descriptionsimplices}, we know that the non-degenerate $m$-simplices of $N\Sigma\cC$ for $m\ge1$ are all of the form
\[x\colon\cO[k]\otimes\cO[\ell]^\circ\to\cC\quad\quad\leftrightarrow\quad\quad x\colon O[k]\otimes O[\ell]^\circ\to C,\]
for some $k,\ell\ge 0$ with $m=k+1+\ell$.
Each such simplex has the following features:
\begin{itemize}[leftmargin=*]
    \item the dimension $m=k+1+\ell$;
    \item the type $k$, which was defined in \cref{definitiontype};
    \item the suspect index $\suspindex$, which will be defined in \cref{DefSuspectIndex}.
\end{itemize}

The goal is to filter the inclusion(s) from \cref{comparisonmap} by a sequence of anodyne extensions $\Sigma\NRS\cC\eqqcolon X_0 \subseteq X_1 \subseteq \ldots \subseteq X_m \subseteq \ldots \subseteq \NRS\Sigma\cC$,
where the inclusions are regular and the underlying simplicial set of $X_m$ contains $X_0$, all the simplices of $\NRS\Sigma\cC$ of dimension less than $m$ as well as all $m+1$-simplices that are \emph{suspect}, a notion that we'll give in \cref{SuspectChar}, and can be computed as
\[k\coloneqq\#\{0\leq i\leq m\ |\ x([i])=\bot\}-1.\]
To this end, we will filter the inclusion $X_{m-1} \subseteq X_{m}$ again via a sequence of anodyne extensions, as $X_{m-1}\eqqcolon  Y_{m-1}\subseteq Y_{m-2} \subseteq \ldots \subseteq Y_{1} \subseteq Y_0=X_m$, where $Y_k$ contains $X_{m-1}$, all the simplices of $\NRS\Sigma\cC$ of dimension $m$ and type at least $k$ as well as all suspect simplices of dimension $m$ and type at least $k-1$. Perhaps surprisingly, to show that the inclusions $Y_{k+1} \subseteq Y_k$ are anodyne, we will construct a further filtration $ Y_{k+1}\eqqcolon W_0\subseteq W_{1} \subseteq \ldots \subseteq W_{k-1} \subseteq W_{k}=Y_k$ and show then that all the inclusions $W_{\suspindex-1}\subseteq W_{\suspindex}$ are anodyne by exhibiting them as a pushout of a sum of anodyne extensions.

\begin{defn}
\label{DefSuspectIndex}
Let $\ell>0$, $k\ge0$ and $m=k+1+\ell$.
Let $x\colon O[k]\otimes O[\ell]^{\circ} \to C$ be a non-totally degenerate
$(k+1+\ell)$-simplex of $N\Sigma\cC$. The \emph{suspect index} of $x$ is the minimal integer $0\leq \suspindex\leq k$ -- if existing -- so that the following two conditions hold.
\begin{enumerate}[label=\normalfont{(SuspInd \arabic*)},leftmargin=*]
\item \label{CondAcross} Whenever $\mathbf{a}\subseteq [\suspindex,k]$, $\mathbf{a'}\subseteq [0,\suspindex-1]$, $\mathbf{b} \subseteq[0,\ell]$, with $\lvert \mathbf{a}\rvert\geq 0$, $\lvert \mathbf{a'}\rvert\geq 0$, $\lvert\mathbf{b}\rvert \geq 0$, we have
    \[
   x([\mathbf{a'},\mathbf{a}] \otimes [0,\mathbf{b}] )=0
    \]
   
    \item\label{CondBelow} Whenever $\textbf a\subseteq[\suspindex, k]$, $\mathbf{b} \subseteq[0,\ell]$, with
 $\lvert \mathbf{b}\vert\geq 1$, we have
    \[
    x([\mathbf{a}]\otimes [\mathbf{b}])=0.
    \]
\end{enumerate}
If there is no integer $\suspindex$ for which the conditions hold, we say that the \emph{suspect index} of $x$ is $k+1$.
\end{defn}

The following are direct consequences of \cite[Proposition 3.4]{SteinerUniversal}.

\begin{lem} 
\label{CriterionDeg}
\label{SuspectChar}
Let $x\colon O[k]\otimes  O[\ell]^{\circ} \to C$ be a non-totally degenerate  $(k+1+\ell)$-simplex of $N\Sigma\cC$.
For $0\leq i\leq k-1$ the following are equivalent:
\begin{enumerate}[leftmargin=*]
    \item Whenever $\mathbf{a'}\subseteq [0,i-1]$, $\mathbf{a}\subseteq [i+2,k]$, with $\lvert \mathbf{a}\rvert\geq -1$, $\lvert \mathbf{a'}\rvert\geq -1$,$\lvert \mathbf{b}\rvert\geq 0$ we have $x([\mathbf{a'},i,i+1,\mathbf{a}] \otimes [\mathbf{b}])=0$.
    \item The $(k+1+\ell)$-simplex $x\colon O[k]\otimes  O[\ell]^{\circ} \to C$ of $(N\Sigma\cC)_{k+1+\ell}$ is degenerate at $i$.
\end{enumerate}
For $k+1\leq i\leq m-1=k+\ell$ the following are equivalent:
\begin{enumerate}[leftmargin=*]
    \item Whenever $\mathbf{a}\subseteq [0,k]$, $\mathbf{b'}\subseteq [0, i-(k+1)-1]$,$\mathbf{b}\subseteq [i-(k+1)+2, \ell]$, with $\lvert \mathbf{a}\rvert\geq 0$, $\lvert \mathbf{b'}\rvert\geq -1$,$\lvert \mathbf{b}\rvert\geq -1$, we have $x([\mathbf{a}]\otimes [\mathbf{b'},i-(k+1),i-(k+1)+1,\mathbf{b}] )=0$.
   \item The $(k+1+\ell)$-simplex $x\colon O[k]\otimes  O[\ell]^{\circ} \to C$ of $(N\Sigma\cC)_{k+1+\ell}$ is degenerate at $i$.
\end{enumerate}
\end{lem}

\begin{lem}
\label{SuspectChar}
Let $y\colon O[k+1]\otimes  O[\ell]^{\circ} \to C$ be a non-totally degenerate  $(k+2+\ell)$-simplex of $N\Sigma\cC$ of suspect index $\suspindex$.
The following are equivalent:
\begin{enumerate}[leftmargin=*]
    \item Whenever $\mathbf{a'}\subseteq [0,\suspindex-2]$, $\mathbf{a}\subseteq [\suspindex+1,m]$, with $\lvert \mathbf{a}\rvert\geq -1$, $\lvert \mathbf{a'}\rvert\geq -1$, we have $y([\mathbf{a'},\suspindex-1,\suspindex,\mathbf{a}] \otimes [0])=0$.
    \item The $(k+1)$-simplex $y([-]\otimes [0])\colon O[k+1] \to \cC$ of $(N\cC)_{k+1}$ is degenerate at $\suspindex-1$.
\end{enumerate}
If the equivalent conditions are satisfied, we say that $y$ is a \emph{suspect} simplex.
\end{lem}

We analyze the faces of a suspect simplex.

\begin{lem}
\label{BoundarySuspect2}
Let $y\colon O[k+1]\otimes  O[\ell]^{\circ} \to C$ be a suspect simplex of $N\Sigma\cC$
suspect index $\suspindex$. Then, if $d_iy$ is not degenerate, we have that
$d_iy$ is 
\begin{enumerate}[leftmargin=*,  label= \normalfont{(Face \arabic*)}]
    \item \label{BeforeMu} 
    of suspect index at most $\suspindex-1$ if $0\leq i \leq \suspindex-2$, 
       \item \label{MuMin1} 
       of suspect index at most $\suspindex-1$ if $i=\suspindex-1$,
          \item \label{AfterMu} 
          either of suspect index at most $\suspindex-1$ or suspect of dimension $m$ and suspect index $\suspindex$ if $\suspindex+1\leq i \leq k+1$,
             \item \label{FacesColumn} 
             of type $k+1$ if $k+2\leq i \leq (k+1)+1+l=m+1$.
\end{enumerate}
\end{lem}
 
\begin{proof}
We distinguish several cases, which correspond to the different cases appearing in the statement.
\begin{enumerate}[leftmargin=0.1cm, label= \normalfont{(Face \arabic*)}, labelwidth=-1.4cm]
    \item  Let $0\leq i \leq \suspindex-2$. Whenever $\mathbf{a'}\subseteq [0,\suspindex-2], \mathbf{a} \subseteq [\suspindex-1,k]$, with
$\lvert \mathbf{a}\rvert\geq 0$, $\lvert \mathbf{a'}\rvert\geq 0$, $\lvert \mathbf{b}\rvert \geq 0$,  
we have $d^i[\mathbf{a'}] \subseteq [0,\suspindex-1]$, $d^i[\mathbf{a}]\subseteq [\suspindex,k+1]$ so that
    \[
    (d_iy)([\mathbf{a'},\mathbf{a}] \otimes [0,\mathbf{b}] )=y(d^i[\mathbf{a'},\mathbf{a}] \otimes [0,\mathbf{b}] )=0,
    \]
    yielding that \ref{CondAcross} holds for $\suspindex-1$ and $d_iy$.
    Whenever $\mathbf{a}\subseteq [\suspindex-1, k]$, with $\lvert \mathbf{a}\vert\geq 0$, $\lvert \mathbf{b}\vert\geq 1$, we have $d^i[\mathbf{a}]\subseteq [\suspindex, k+1]$ and so
    \[
    (d_iy)([\mathbf{a}]\otimes [\mathbf{b}])=y(d^i[\mathbf{a}]\otimes [\mathbf{b}])=0,
    \]
    yielding that \ref{CondBelow} holds for $\suspindex-1$ and $d_iy$.
    So the suspect index of $d_iy$ is at most $\suspindex-1$.

\item Let $i=\suspindex-1$. Whenever $\mathbf{a'}\subseteq [0,\suspindex-2], \mathbf{a} \subseteq [\suspindex-1,k]$, with $\lvert \mathbf{a}\rvert\geq 0$, $\lvert \mathbf{a'}\rvert\geq 0$, $\lvert \mathbf{b}\rvert \geq 0$, 
we have $d^{\suspindex-1}[\mathbf{a'}] \subseteq [0,\suspindex-1]$, $d^{\suspindex-1}[\mathbf{a}]\subseteq [\suspindex,k+1]$
so that
    \[
    (d_iy)([\mathbf{a'},\mathbf{a}] \otimes [0,\mathbf{b}] )=y(d^i[\mathbf{a'},\mathbf{a}] \otimes [0,\mathbf{b}] )=0,
    \]
     yielding that \ref{CondAcross} holds for $\suspindex-1$ and $d_iy$.
Whenever $\mathbf{a}\subseteq [\suspindex-1, k]$ , with
$\lvert \mathbf{a}\vert\geq 0$, $\lvert \mathbf{b}\vert\geq 1$, we have $d^{\suspindex-1}\mathbf{a}\subseteq [\suspindex, k+1]$ and so
    \[
    (d_{i}y)([\mathbf{a}]\otimes [\mathbf{b}])=y(d^i[\mathbf{a}]\otimes [\mathbf{b}])=0,
    \]
yielding that \ref{CondBelow} holds for $\suspindex-1$ and $d_iy$. So the suspect index of $d_iy$ is at most $\suspindex-1$.

\item Let $\suspindex+1\leq i \leq k+1$. Whenever $\mathbf{a'}\subseteq [0,\suspindex-1], \mathbf{a} \subseteq [\suspindex,k]$, with
$\lvert \mathbf{a}\rvert\geq 0$, $\lvert \mathbf{a'}\rvert\geq 0$, $\lvert \mathbf{b}\rvert \geq 0$,  
we have
     \[
    (d_iy)([\mathbf{a}]\otimes [\mathbf{b}])=0,
    \]
    yielding that \ref{CondAcross} holds for $\suspindex$ and $d_iy$.
    
    Moreover, whenever
$\lvert \mathbf{a}\vert\geq 0$, $\mathbf{a} \subseteq [\suspindex,k]$, with $\lvert \mathbf{b}\vert\geq 1$, we have
    \[
    (d_{i}y)([\mathbf{a}]\otimes [\mathbf{b}])=0,
    \]
    yielding that \ref{CondBelow} holds for $\suspindex$ and $d_iy$. So the suspect index of $d_iy$ is at most $\suspindex$.

To see that the simplex is suspect in the case that the suspect index is exactly $\suspindex$, we observe that for any
    $\mathbf{a'}\subseteq [0,\suspindex-2], \mathbf{a} \subseteq [\suspindex-1,k]$,
    with $\lvert \mathbf{a}\rvert\geq -1$, $\lvert \mathbf{a'}\rvert\geq -1$ we have
    \[
     (d_iy)([\mathbf{a'},\suspindex-1,\suspindex,\mathbf{a}]\otimes [0])= y(d^i[\mathbf{a'},\suspindex-1,\suspindex,\mathbf{a}] \otimes [0])=0.
    \]
\item The last case is clear since the face operator acts on the second coordinate. 
\qedhere
\end{enumerate}
\end{proof}

\begin{lem}
\label{parent}
Let $x$ be a non-degenerate non-suspect simplex of suspect index $\suspindex$. There is a simplex $\widetilde x\colon O[k+1]\otimes O[\ell]^{\circ}\to C$ defined by the following formulas:
\begin{enumerate}[leftmargin=*, label= \normalfont{(P\arabic*)}]
\def\boxlength{7.9cm}
    \item \label{SuspFaceCorrect} \makebox[\boxlength][l]{$\widetilde x(d^{\suspindex}[\mathbf{a}]\otimes[\mathbf{b}]) = x([\mathbf{a}]\otimes[\mathbf{b}])$} if $\lvert \mathbf{a}\rvert\geq -1, \lvert \mathbf{b}\rvert\geq -1$, 
   \item \label{SuspColumnBelow}  \makebox[\boxlength][l]{$\widetilde x([\suspindex,\mathbf{a}]\otimes [b] )=x(s^{\suspindex-1}[\suspindex,\mathbf{a}]\otimes [0])$}\mbox{ if }$ \lvert\mathbf{a}\rvert \geq -1$;
   \item \label{SuspArrayAbove}\makebox[\boxlength][l]{$\widetilde x([\mathbf{a'},\suspindex] \otimes [\mathbf{b}])= x([\mathbf{a'}] \otimes [0,\mathbf{b}])$}  if $\lvert \mathbf{a'}\vert\geq 0,\lvert \mathbf{b}\rvert \geq 1$;
   \item \label{SuspColumnAbove}\makebox[\boxlength][l]{$\widetilde x([\mathbf{a'},\suspindex] \otimes [b])= x([\mathbf{a'}] \otimes [0,b])+x([\mathbf{a'},\suspindex-1] \otimes [0])$} if $\lvert \mathbf{a'}\vert\geq 0$;
   \item \label{SuspColumnAcross}\makebox[\boxlength][l]{$\widetilde x([\mathbf{a'},\suspindex,\mathbf{a}] \otimes [b]) =x(s^{\suspindex-1}[\mathbf{a'},\suspindex,\mathbf{a}] \otimes [0])$} if $\lvert \mathbf{a}\rvert\geq 0, \lvert \mathbf{a'}\rvert\geq 0$;
   \item \label{SuspRowAcross}\makebox[\boxlength][l]{$\widetilde x[\mathbf{a'},\suspindex,\mathbf{a}] \otimes [\mathbf{b}] = 0$} if $\lvert \mathbf{a}\rvert\geq 0, \lvert \mathbf{a'}\rvert\geq -1, \lvert \mathbf{b}\rvert \geq 1$;
   \item \label{SuspCritRow}\makebox[\boxlength][l]{$\widetilde x([\suspindex] \otimes [\mathbf{b}] )=0$} if $ \lvert \mathbf{b}\rvert \geq 1$.
\end{enumerate}

\end{lem}

The proof consists of a careful (and tedious) analysis of all possible cases, and is postponed until \cref{appendix}.

\begin{rmk}
\label{facextilde}
Given a non-degenerate non-suspect simplex $x$ of suspect index $\suspindex$, by construction we have $d_\suspindex(\widetilde x)=x$.
\end{rmk}

We record the following features of $\widetilde x$.

\begin{lem}
\label{ParametersTilde}
If $x$ is a non-suspect simplex of $N\Sigma\cC$ with dimension $m=k+1+\ell$, type $k$ and suspect index $\suspindex$, then the simplex $\widetilde x$ is a suspect simplex of dimension $m+1$, type $k+1$ and suspect index $\suspindex$.
\end{lem}

\begin{proof}
By construction, the simplex $\widetilde x$ is a suspect simplex of dimension $m$, type $k$ and suspect index at most $\suspindex$. We now prove the suspect index of $\widetilde x$ is exactly $\suspindex$.

Assume by contradiction that the suspect index of $\widetilde x$ is at most $\suspindex-1$. Then, whenever $\mathbf{a'} \subseteq [0,\suspindex-2]$, $\mathbf{a} \subseteq [\suspindex-1,k]$, we have $d^{\suspindex}[\mathbf{a}] \subseteq [\suspindex-1,k]$, so that 
\[
x([\mathbf{a'},\mathbf{a}] \otimes [0,\mathbf{b}] )=\widetilde x(d^{\suspindex}[\mathbf{a'},\mathbf{a}] \otimes [0,\mathbf{b}] ))=0,
\]
 yielding that \ref{CondAcross} holds for $\suspindex-1$ and $x$. 
Also, whenever $\mathbf{a}\subseteq[\suspindex-1,k]$, with $\lvert \mathbf{a}\vert\geq 0$,  
    and $\lvert \mathbf{b}\vert\geq 1$, we have
    \[
    x([\mathbf{a}]\otimes [\mathbf{b}])=\widetilde x(d^{\suspindex}[\mathbf{a}]\otimes [\mathbf{b}])=0,
    \]
     yielding that \ref{CondBelow} holds for $\suspindex-1$ and $x$. This would thus imply that $x$ is also of suspect index at most $\suspindex-1$, contrary to the assumption.
\end{proof}

\begin{lem}
\label{xtildenondegenerate}
Let $x$ be a non-suspect simplex of $N\Sigma\cC$. If $x$ is non-degenerate, then $\widetilde x$ is non-degenerate.
\end{lem}

\begin{proof}
Let $x$ be a simplex of dimension $m$, type $k$ and suspect index $\suspindex$. Assume that
$\widetilde x= s_i z$
is degenerate at some $0\leq i\leq m$, and deduce a contradiction by distinguishing several cases.
\begin{itemize}[leftmargin=*]

    \item If $i=\suspindex-1$, we can prove that $x$ would be of suspect index at most $\suspindex-1$, contradicting 
    the assumption.
    Since $\widetilde x$ is degenerate at $\suspindex-1$, we have $x=d_{\suspindex}\widetilde x=d_{\suspindex-1}\widetilde x$.
    
    We check first that then the condition \ref{CondAcross} for $\suspindex-1$ and $x$. Assume $\mathbf{a'}\subseteq [0,\suspindex-2], \mathbf{a} \subseteq [\suspindex-1,k]$, with $\lvert \mathbf{a}\rvert\geq 0$, $\lvert \mathbf{a'}\rvert\geq 0$, $\lvert \mathbf{b}\rvert \geq 0$. If $\mathbf{a}$ does not contain $\suspindex-1$, then we have 
    \[
  x([\mathbf{a'},\mathbf{a}] \otimes [0,\mathbf{b}] )=0
    \]
since $x$ has suspect index $\suspindex$. If $\mathbf{a}$ contains $\suspindex-1$ and at least one further element, then
we have 
    \[
  x([\mathbf{a'},\suspindex,\mathbf{a}] \otimes [0,\mathbf{b}] )=0
    \]
since $x$ has suspect index $\suspindex>\suspindex-1$.
Finally, if $\mathbf{a}=[\suspindex-1]$, we have
     \[
     \begin{array}{llllll}
  x([\mathbf{a'},\suspindex-1] \otimes [0,\mathbf{b}] )&=&d_{\suspindex-1}\widetilde x([\mathbf{a'},\suspindex-1] \otimes [0,\mathbf{b}] )\\
  &=&\widetilde x([\mathbf{a'},\suspindex] \otimes [0,\mathbf{b}] )=0.
  \end{array}
    \]
    This yields the condition \ref{CondAcross} for $\suspindex-1$ and $x$. 
    
For \ref{CondBelow}, we observe that, whenever $\mathbf{a} \subseteq [\suspindex,k]$, we have
\[x([\suspindex-1, \mathbf{a}] \otimes [\mathbf{b}])=\widetilde x([\suspindex,\mathbf{a}] \otimes [\mathbf{b}])=0.\]
    This yields the condition \ref{CondBelow} for $\suspindex-1$ and $x$. 
    \item  If $i=\suspindex$, we can prove that $x$ would be suspect, contradicting the assumption. Indeed, we observe that
         \[
     \begin{array}{llllll}
  x([\mathbf{a'},\suspindex-1, \suspindex, \mathbf{a}] \otimes [0] )&=& x(s^{\suspindex-1}d^{\suspindex-1}([\mathbf{a'},\suspindex-1, \suspindex, \mathbf{a}] \otimes [0] )& \\
  &=& \widetilde{x}(d^{\suspindex-1}([\mathbf{a'},\suspindex-1, \suspindex, \mathbf{a}] \otimes [0] ) &\text{\ref{SuspColumnAcross}, \ref{SuspColumnBelow}}\\  &=& (s_{\suspindex}z)(d^{\suspindex-1}([\mathbf{a'},\suspindex-1, \suspindex, \mathbf{a}] \otimes [0] )&\text{Assumption}\\
  &=&0&\text{\cref{SuspectChar}}
  \end{array}
    \]
    so that $x$ would be indeed suspect. This yields the desired contradiction.

        \item  If $i\neq\suspindex-1,\suspindex$, we can prove that $x$ would be degenerate, contradicting the assumption. Indeed, we would obtain
 \[
 x=d_{\suspindex} \widetilde x = d_{\suspindex} s_i z=\left\{
 \begin{array}{llll}
  s_{i-1} d_{\suspindex} z&\text{if }\suspindex+1\leq i\leq k+1+l\\
 s_{i} d_{\suspindex} z, &\text{if }0\leq i\leq \suspindex-2
 \end{array}
 \right.
 \]
 contrary to the assumption that $x$ is non-degenerate.
\end{itemize}
This concludes the proof.
\end{proof}

The following shows that the values of a suspect simplex $y$ are enforced by those of its face $d_{\suspindex}y$.

\begin{lem}
\label{PartnerEnforced}
Let $y\colon O[k+1]\otimes  O[\ell]^{\circ} \to C$ be a suspect simplex of $N\Sigma\cC$
with suspect index $\suspindex$.
\begin{enumerate}[leftmargin=*, label= \normalfont{(S\arabic*)}]
\def\boxlength{8.3cm}
    \item \label{P1} \makebox[\boxlength][l]{$y([d^{\suspindex}\mathbf{a}]\otimes[\mathbf{b}]) = d_{\suspindex}y([\mathbf{a}]\otimes[\mathbf{b}])$} if $\lvert \mathbf{a}\rvert\geq -1, \lvert \mathbf{b}\rvert\geq -1$, 
   \item  \label{P2}\makebox[\boxlength][l]{$y([\suspindex,\mathbf{a}]\otimes [b] )=d_{\suspindex}y(s^{\suspindex-1}[\suspindex,\mathbf{a}]\otimes [0])$}\mbox{ if }$ \lvert\mathbf{a}\rvert \geq -1$;
   \item \label{P3}\makebox[\boxlength][l]{$y([\mathbf{a'},\suspindex] \otimes [\mathbf{b}])= d_{\suspindex}y([\mathbf{a'}] \otimes [0,\mathbf{b}])$}  if $\lvert \mathbf{a'}\vert\geq 0,\lvert \mathbf{b}\rvert \geq 1$;
   \item \label{P4}\makebox[\boxlength][l]{$y([\mathbf{a'},\suspindex] \otimes [b])= d_{\suspindex}y([\mathbf{a'}] \otimes [0,b])+d_{\suspindex}y([\mathbf{a'},\suspindex-1] \otimes [0])$} if $\lvert \mathbf{a'}\vert\geq 0$;
   \item \label{P5}\makebox[\boxlength][l]{$y([\mathbf{a'},\suspindex,\mathbf{a}] \otimes [b]) =d_{\suspindex}y(s^{\suspindex-1}[\mathbf{a'},\suspindex,\mathbf{a}] \otimes [0])$} if $\lvert \mathbf{a}\rvert\geq 0, \lvert \mathbf{a'}\rvert\geq 0$;
   \item \label{P6}\makebox[\boxlength][l]{$y([\mathbf{a'},\suspindex,\mathbf{a}] \otimes [\mathbf{b}]) = 0$} if $\lvert \mathbf{a}\rvert\geq 0, \lvert \mathbf{a'}\rvert\geq -1, \lvert \mathbf{b}\rvert \geq 1$;
   \item \label{P7} \makebox[\boxlength][l]{$y([\suspindex] \otimes [\mathbf{b}] )=0$} if $ \lvert \mathbf{b}\rvert \geq 1$.
\end{enumerate}

\end{lem}

\begin{proof}
We now show that a suspect $m$-simplex $y\colon O[k+1]\otimes O[\ell]^{\circ} \to C$ of suspect index $\suspindex$ is already completely specified by $d_{\suspindex} y$, in the way described more precisely by the statement. 
\begin{enumerate}[leftmargin=0.1cm, label=(S\arabic*), labelwidth=-0.8cm]
    \item The value of $y$ in this case is by definition of simplicial structure of $N\Sigma\cC$.
    \item We prove the formula for $y([\suspindex,\mathbf{a}]\otimes [b])$ in this case.
    We observe that 
\[
\begin{array}{llll}
0&=&\partial y([\suspindex,\mathbf{a}]\otimes [0,b])& \text{\ref{CondBelow}}\\
&=& y([\mathbf{a}]\otimes [0,b])-y([\suspindex,\partial\mathbf{a}]\otimes [0,b])&\\
&&+(-1)^{\lvert \mathbf{a}\rvert+1}
y([\suspindex,\mathbf{a}]\otimes [0])-(-1)^{\lvert \mathbf{a}\rvert+1}y([\suspindex,\mathbf{a}]\otimes [b])&\\
&=&(-1)^{\lvert \mathbf{a}\rvert+1}
y([\suspindex,\mathbf{a}]\otimes [0])-(-1)^{\lvert \mathbf{a}\rvert+1}y([\suspindex,\mathbf{a}]\otimes [b]) &\text{\ref{CondBelow}}\\
&=&y([\suspindex,\mathbf{a}]\otimes [0])-y([\suspindex,\mathbf{a}]\otimes [b])&\\
&=&y(d^{\suspindex}s^{\suspindex-1}[\suspindex,\mathbf{a}]\otimes [0])-y([\suspindex,\mathbf{a}]\otimes [b])&\text{\cref{SuspectChar}}
\end{array}
\]
The desired formula follows.

    \item We prove the formula for $y([\mathbf{a'}, \suspindex]\otimes [\mathbf{b}])$ in this case. If $\mathbf{b}$ contains $0$, the term vanishes and the formula follows by \ref{CondAcross}. If $\mathbf{b}$ does not contain $0$, we have
    \[
\begin{array}{llll}
0&=&\partial y([\mathbf{a'},\suspindex]\otimes [0,\mathbf{b}])& \text{\ref{CondAcross}}\\
&=& (-1)^{\lvert \mathbf{a'}\rvert+1}y([\mathbf{a'}]\otimes [0,\mathbf{b}])+y([\partial\mathbf{a'},\suspindex]\otimes [0,\mathbf{b}])&\\
&&+(-1)^{\lvert \mathbf{a'}\rvert+2}y([\mathbf{a'},\suspindex]\otimes [\mathbf{b}])
+(-1)^{\lvert \mathbf{a'}\rvert+2}y([\mathbf{a'},\suspindex]\otimes [0,\partial^{\circ} \mathbf{b}])&\\
&=&(-1)^{\lvert \mathbf{a'}\rvert+1}y([\mathbf{a'}]\otimes [0,\mathbf{b}])+(-1)^{\lvert \mathbf{a'}\rvert+2}y([\mathbf{a'},\suspindex]\otimes [\mathbf{b}])&\text{\ref{CondAcross},\ref{CondBelow}}\\
&=&y([\mathbf{a'}]\otimes [0,\mathbf{b}])-y([\mathbf{a'},\suspindex]\otimes [\mathbf{b}])&
\end{array}
\]
The desired formula follows. 
    \item We prove the formula for $y([\mathbf{a'},\suspindex] \otimes [b])$
    in this case. If $b=0$, the formula follows from \cref{SuspectChar}.
    If $b>0$, 
 we have
\[
\begin{array}{llll}
0&=&\partial y([\mathbf{a'},\suspindex]\otimes [0,b])&\text{\ref{CondAcross}}\\
&=& (-1)^{\lvert \mathbf{a'}\rvert+1}y([\mathbf{a'}]\otimes [0,b])
+y([\partial\mathbf{a'},\suspindex]\otimes [0,b])&\\
&&+(-1)^{\lvert \mathbf{a'}\rvert+2}y([\mathbf{a'},\suspindex]\otimes [b])
+(-1)^{\lvert \mathbf{a'}\rvert+3}y([\mathbf{a'},\suspindex]\otimes [0])&\\
&=& (-1)^{\lvert \mathbf{a'}\rvert+1}y([\mathbf{a'}]\otimes [0,b])\\
&&+(-1)^{\lvert \mathbf{a'}\rvert+2}y([\mathbf{a'},\suspindex]\otimes [b])
+(-1)^{\lvert \mathbf{a'}\rvert+3}y([\mathbf{a'},\suspindex]\otimes [0])&\text{\ref{CondAcross},\ref{CondBelow}}\\
&=& y/[\mathbf{a'}]\otimes [0,b])
-y([\mathbf{a'},\suspindex]\otimes [b])
+y([\mathbf{a'},\suspindex]\otimes [0])&\\
&=& y([\mathbf{a'}]\otimes [0,b])
-y([\mathbf{a'},\suspindex]\otimes [b])
+y([\mathbf{a'},\suspindex-1]\otimes [0])&\text{\cref{SuspectChar}}.
\end{array}
\]
The desired formula follows. 
    \item We prove the formula for $y([\mathbf{a'},\suspindex,\mathbf{a}] \otimes [b])$ in this case.
    If $b=0$, then the equality follows from \cref{SuspectChar}. If $b>0$, we have
\[
\begin{array}{llll}
0&=&\partial y([\mathbf{a'},\suspindex,\mathbf{a}] \otimes [0,b])&\text{\ref{CondAcross}}\\
&=&y([\partial\mathbf{a'},\suspindex,\mathbf{a}] \otimes [0,b]) + (-1)^{\lvert \mathbf{a'}\rvert+1}y([\mathbf{a'},\mathbf{a}] \otimes [0,b])& \\
&&+(-1)^{\lvert \mathbf{a'}\rvert+2}y([\mathbf{a'},\suspindex,\partial\mathbf{a}] \otimes [0,b])+(-1)^{\lvert \mathbf{a'}\rvert+1+\lvert \mathbf{a}\rvert}y([\mathbf{a'},\suspindex,\mathbf{a}] \otimes [0])&\\
&&-(-1)^{\lvert \mathbf{a'}\rvert+1+\lvert \mathbf{a}\rvert}y([\mathbf{a'},\suspindex,\mathbf{a}] \otimes [b])&\\
&=&(-1)^{\lvert \mathbf{a'}\rvert+1+\lvert \mathbf{a}\rvert}y([\mathbf{a'},\suspindex,\mathbf{a}] \otimes [0])-(-1)^{\lvert \mathbf{a'}\rvert+1+\lvert \mathbf{a}\rvert}y([\mathbf{a'},\suspindex,\mathbf{a}] \otimes [b])&\text{\ref{CondAcross},\ref{CondBelow}}\\
&=&y([\mathbf{a'},\suspindex,\mathbf{a}] \otimes [0])-y([\mathbf{a'},\suspindex,\mathbf{a}] \otimes [b])&\\
&=&y([\mathbf{a'},\suspindex,\mathbf{a}] \otimes [0])-d_{\suspindex}y(s^{\suspindex-1}[\mathbf{a'},\suspindex,\mathbf{a}] \otimes [0])&\text{\cref{SuspectChar}}.
\end{array}
\]
The desired formula follows. 
\item We prove that $y([\mathbf{a'},\suspindex,\mathbf{a}] \otimes [\mathbf{b}])$ necessarily vanishes. If $\mathbf{b}$ contains $0$, this is a special case of \ref{CondAcross}. If $\mathbf{b}$ does not contain $0$, we have
    \[
\begin{array}{llll}
0&=&\partial y([\mathbf{a'},\suspindex,\mathbf{a}]\otimes [0,\mathbf{b}])& \text{\ref{CondAcross}}\\
&=& y([\partial\mathbf{a'},\suspindex,\mathbf{a}]\otimes [0,\mathbf{b}])
+(-1)^{\lvert \mathbf{a'}\rvert+1}y([\mathbf{a'},\mathbf{a}]\otimes [0,\mathbf{b}])&\\
&&+(-1)^{\lvert \mathbf{a'}\rvert+2}y([\mathbf{a'},\suspindex,\partial\mathbf{a}]\otimes [0,\mathbf{b}])&\\
&&+(-1)^{\lvert \mathbf{a'}\rvert+2+\lvert \mathbf{a}\rvert}y([\mathbf{a'},\suspindex,\mathbf{a}]\otimes [\mathbf{b}])
+(-1)^{\lvert \mathbf{a'}\rvert+2+\lvert \mathbf{a}\rvert}y([\mathbf{a'},\suspindex,\mathbf{a}]\otimes [0,\partial^{\circ} \mathbf{b}])&\\
&=&(-1)^{\lvert \mathbf{a'}\rvert+2}y([\mathbf{a'},\suspindex,\partial\mathbf{a}]\otimes [0,\mathbf{b}])+(-1)^{\lvert \mathbf{a'}\rvert+2+\lvert \mathbf{a}\rvert}y([\mathbf{a'},\suspindex,\mathbf{a}]\otimes [\mathbf{b}])&\text{\ref{CondAcross},\ref{CondBelow}}
\end{array}
\]
If $\lvert \mathbf{a}\rvert>0$, the desired vanishing follows directly from \ref{CondAcross}. If $\lvert \mathbf{a}\rvert=0$, the desired vanishing follows from \ref{P3} which we have treated before.
\item The fact that $y([\suspindex] \otimes [\mathbf{b}] )$ vanishes in this case can be seen applying \ref{CondBelow} for $y$ and $\suspindex$.
\end{enumerate}

This concludes the proof.
\end{proof}

\begin{lem}
\label{featuresface}
Let $y\colon O[k+1]\otimes  O[\ell]^{\circ} \to C$ be a non-degenerate suspect simplex of $N\Sigma\cC$
with suspect index $\suspindex$.
Then the face $d_{\suspindex}y$ is a simplex of dimension $m$, type $k$ and suspect index $\suspindex$.
\end{lem}

We record the following features of $d_{\suspindex}y$.

\begin{proof}
The value of the dimension and type of $d_\suspindex y$ are immediate from the definitions.

We now show that $d_{\suspindex}y$ is of suspect index $\suspindex$. 
It is straightforward from the construction that the suspect index of $d_\suspindex y$ is at most $\suspindex$. Assuming for contradiction the suspect index of $d_{\suspindex}y$ to be at most $\suspindex-1$, we show that the suspect index of $y$ would be also at most $\suspindex-1$, contrary to the assumptions.

One can show that the conditions \ref{CondAcross} and \ref{CondBelow} hold for $y$ and $\suspindex-1$. If $[\mathbf a]$ contains $\suspindex$, this follows from \ref{CondAcross} and \ref{CondBelow}. If $[\mathbf a]$ does not contain $\suspindex$, this follows from \ref{P3}, \ref{P6}, \ref{P7}
of \cref{PartnerEnforced}.

Finally, we show that $d_{\suspindex}y$ is a non-suspect simplex. Assuming for contradiction that $d_{\suspindex}y$ is suspect, we show using the characterization from \cref{CriterionDeg} that then $y$ is degenerate at $\suspindex$, contrary to the assumptions.
To this end, we need to show that $y([\mathbf{a'}, \suspindex, \suspindex +1, \mathbf{a}]\otimes [\mathbf{b}])$ vanishes. If $\lvert \mathbf{b}\rvert\geq 1$, this term vanishes by \ref{P6}. If $\lvert \mathbf{b}\rvert=0$ and $\lvert \mathbf{a'}\rvert=-1$, then this term vanishes by \ref{P2} using the assumption that $d_{\suspindex}y$ is suspect with suspect index $\suspindex$. If $\lvert \mathbf{b}\rvert=0$ and $\lvert \mathbf{a'}\rvert\geq 0$, this term vanishes by \ref{P5} using the assumption that $d_{\suspindex}y$ is suspect.
\end{proof}

\begin{lem}
\label{facenondegenerate}
Let $y$ be a suspect simplex of $N\Sigma\cC$ with suspect index $\suspindex$. If $y$ is non-degenerate, then the face $d_{\suspindex}y$ is a non-degenerate simplex.
\end{lem}

\begin{proof}
Let $k+1$ be the type of $y$. Assuming for contradiction that $d_{\suspindex}y$ is degenerate, we show that $y$ itself has to be degenerate at some $0\leq i\leq k+\ell$. Notice that the case $i=k$ cannot happen, because it would violate the type property from \cref{featuresface}.
\begin{itemize}[leftmargin=*]
    \item If $0\leq i <\suspindex-1$, we use \cref{CriterionDeg} to show that $y$ is degenerate at $i$. The fact that
    \[y([\mathbf{a'},i,i+1,\mathbf{a}]\otimes [\mathbf{b}])\]
    vanishes is by definition when $\suspindex$ does not occur in $[\mathbf a]$. Otherwise, it can be deduced from the formulas \ref{P3}, \ref{P4}, \ref{P5}, \ref{P6}, together with the assumption that $d_\suspindex y$ is degenerate at $i$.
\item If $i=\suspindex-1$, we use \cref{CriterionDeg} to show that $y$ is degenerate at $\suspindex$. The fact that
\[
y([\mathbf{a'},\suspindex,\suspindex+1,\mathbf{a}]\otimes [\mathbf{b}])
\]
vanishes can be deduced from the formulas \ref{P2}, \ref{P5}, \ref{P6}, together with the assumption that $d_\suspindex y$ is degenerate at $\suspindex-1$ in the formulation from \cref{CriterionDeg}.

    \item  If $\suspindex\leq i <  k$, we use \cref{CriterionDeg} to show that $y$ is degenerate  at $i+1$. The fact that
\[
y([\mathbf{a'},i+1,i+2,\mathbf{a}]\otimes [\mathbf{b}]).
\]
vanishes can be deduced from the formulas \ref{P2}, \ref{P5}, \ref{P6}, together with the assumption that $d_\suspindex y$ is degenerate at $i$ in the formulation from \cref{CriterionDeg}.
\item  If $k+1\leq i \leq k+\ell$, we use \cref{CriterionDeg} to show that $y$ is degenerate at $i+1$. The fact that
\[y([\mathbf a]\otimes[\mathbf b])\]
vanishes follows from the formulas from \cref{PartnerEnforced}, together with the fact that $d_\suspindex y$ is degenerate at $i$ in the formulation from \cref{CriterionDeg}.
\end{itemize}
This concludes the proof.
\end{proof}


We can now establish a correspondence between the suspect and non-suspect simplices of $N\Sigma\cC$.

\begin{prop}
\label{lemmabijection}
Let $\cC$ be a $1$-category and $m\ge 0$. Recall the inclusion $\Sigma N\cC\hookrightarrow  N\Sigma\cC$ from \cref{comparisonmap}.
\begin{enumerate}[label=(\roman*),leftmargin=0pt, align=left]
    \item\label{lemmabijectionI} The non-degenerate $(m+1)$-simplices in $\Sigma N\cC$, regarded as a simplicial subset of $N\Sigma\cC$, are contained in the $(m+1)$-simplices of type $m$.
    \item The non-degenerate $(m+1)$-simplices in $N\Sigma\cC$ that do not belong to $\Sigma N\cC$ are non-degenerate $(m+1)$-simplices of type $0\leq k\leq m-1$ and suspect index $1\leq \suspindex\leq k+1$. 
    \item The $\suspindex$-th face map gives a bijective correspondence between
    the non-degenerate suspect $(m+1)$-simplices $\widetilde{x}$ in $\Sigma N\cC\setminus \Sigma N\cC$ of type $1\leq k\leq m-1$ and suspect index $1 \leq \suspindex \leq k+1$ and the non-degenerate non-suspect $m$-simplices $x$ of type $0\leq k-1\leq m-2$ and suspect index $1\leq \suspindex\leq k+1$.
\end{enumerate}
\end{prop}

\begin{proof}
The first two statements can be verified by direct inspection. For the third statement, we observe that the assignment $\widetilde{(-)}$ from \cref{parent} is an inverse for the function $d_{\suspindex}$ with the given domain and codomain. Indeed, \cref{featuresface,facenondegenerate} show that $y\mapsto d_{\suspindex}y$ is a well-defined function, \cref{ParametersTilde,xtildenondegenerate} show that $x\mapsto\widetilde x$ is a well-defined function, \cref{facextilde} shows that $d_\suspindex\widetilde x=x$, and the formulas from \cref{parent,PartnerEnforced} together imply that $\widetilde{d_\suspindex y}=y$.
\end{proof}

We now prove the theorem.

\begin{proof}[Proof of \cref{ThmB}]
We prove (2), and (1) follows then by applying the forgetful functor from marked simplicial sets to simplicial sets. In order to show that the inclusion $\Sigma\NRS\cC\to \NRS\Sigma\cC$ is a complicial inner anodyne extension, we will realize it as a transfinite composite of intermediate complicial inner anodyne extensions
\[\Sigma\NRS\cC\eqqcolon X_1\hookrightarrow X_2\hookrightarrow\dots\hookrightarrow X_{m-1}\hookrightarrow X_{m}\hookrightarrow\dots \hookrightarrow   \NRS\Sigma\cC.\]
For $m\geq 2$, we let $X_m$ be the smallest regular subsimplicial set of $\NRS\Sigma\cC$ 
containing $X_{m-1}$, all $m$-simplices of $N\Sigma\cC$, as well as
the suspect $(m+1)$-simplices of $N\Sigma\cC$.
Note that $X_1$ already contains all non-degenerate $1$-simplices of $N\Sigma\cC$ and that there are no non-degenerate suspect $2$-simplices. Moreover, by \cref{lemmabijection}, the subsimplicial set $X_1$ contains all non-degenerate $(m+1)$-simplices of type $m$. 
We see that the difference between $X_{m-1}$ and $X_m$ are the non-degenerate non-suspect $m$-simplices of type at most $m-2$ and the non-degenerate suspect $(m+1)$-simplices of type at most $m-1$.

In order to show that the inclusion $X_{m-1}\hookrightarrow X_{m}$ is a complicial inner anodyne extension for all $d\geq 2$, we realize it as a composite of intermediate complicial inner anodyne extensions
\[X_{m-1}\eqqcolon Y_{m} \hookrightarrow Y_{m-1} \hookrightarrow \ldots \hookrightarrow Y_{k+1}\hookrightarrow Y_{k}\hookrightarrow\ldots \hookrightarrow Y_{1}=X_{m}.\]
For $1\le k< m$,
let $Y_k$ be the
smallest regular subset of $X_{m}$ containing $Y_{k+1}$ as well as 
all non-degenerate suspect $(m+1)$-simplices $\widetilde{ x}$ of $N\Sigma\cC$ of type $k$ and all non-degenerate non-suspect $m$-simplices of type $k-1$. Note that $Y_m$ already contains all non-degenerate $m$-simplices of type $m-1$ and that any suspect $(m+1)$-simplex of type $m$ is necessarily a degeneracy of a $m$-simplex of type $m-1$ and thus can be checked to be also already in $Y_m$.
We see using \cref{BoundarySuspect2,lemmabijection} that the difference between $Y_{k}$ and $Y_{k+1}$ are
the non-degenerate suspect $(m+1)$-simplices of type $k$ and possibly some of their faces (precisely those that are neither suspect nor of type $k$ or higher).

In order to show that the inclusion $Y_{k+1}\hookrightarrow Y_k$ is a complicial inner anodyne extension for $1\le k\le m-1$, we realize it as a filtration made by intermediate complicial inner anodyne extensions
\[Y_{k+1}\eqqcolon W_0\hookrightarrow W_1 \hookrightarrow \ldots \hookrightarrow W_{\suspindex-1}\hookrightarrow W_{\suspindex}\hookrightarrow\ldots\hookrightarrow W_{k}=Y_{k}.\]
For $0<\suspindex\leq k$, we let $W_{\suspindex}$ be the smallest regular simplicial subset of $Y_k$
containing $W_{\suspindex-1}$ as well as 
all suspect $(m+1)$-simplices of $  \NRS\Sigma\cC$ of type $k$ and suspect index $\suspindex$, namely those $\widetilde{ x}$
for which each $i$-th row constant for $\suspindex\leq i \leq k$.
Note that any $m$-simplex of suspect index $0$ is either degenerate or of type $m-1$ and thus can be checked to be already in $X_1\subseteq W_0$.
We see using \cref{BoundarySuspect2,lemmabijection} that the difference between $W_{\suspindex-1}$ and $W_\suspindex$ are the non-degenerate suspect $(m+1)$-simplices $\widetilde{ x}$ of type $k$ and suspect index $\suspindex$ and the non-degenerate non-suspect $m$-simplices $ x$
of type $k-1$ and suspect index $\suspindex$. 
There is a bijective correspondence between the $(m+1)$- and $m$-simplices mentioned above, as shown in \cref{lemmabijection}.
 
We now record some relevant properties of the $(m+1)$-simplices $\widetilde x$ as above. 
\begin{itemize}[leftmargin=*]
    \item We argue by induction and using \cref{BoundarySuspect2} that the $\suspindex$-horn of $\widetilde{ x}$ belongs to $W_{\suspindex-1}$; in particular, the $\suspindex$-horn defines a map of (underlying) simplicial sets
    \[\Lambda^{\suspindex}[m+1]\to W_{\suspindex-1}.\]
Indeed, using \cref{BoundarySuspect2} we see that the $a$-th face of $\widetilde{ x}$ is already in $W_{\suspindex-1}$ for $a\neq \suspindex$ since it is either a degeneracy of a simplex of smaller dimension or: 
\begin{itemize}[leftmargin=*]
\renewcommand\labelitemii{$\diamondsuit$}
    \item if $0\leq a \leq \suspindex-2$, the face $ d_a\widetilde x$ is of suspect index at most $\suspindex-1$, and in particular it belongs to $W_{\suspindex-1}$.
     \item if $a=\suspindex-1$, the face $ d_a\widetilde x$ has suspect index at most $\suspindex-1$, and in particular it belongs to $W_{\suspindex-1}$.
        \item if $\suspindex+1\leq a \leq k$, the face $ d_a\widetilde x$ is either of suspect index at most $\suspindex-1$ or suspect of dimension $m$ and suspect index $\suspindex$; in either case, it belongs to $W_{\suspindex-1}$.%
     \item if $k+1\leq a\leq m+1$, the face $ d_a\widetilde x$ is of type $k+1$,
     and in particular it belongs to $Y_{k+1}\subseteq W_{\suspindex-1}$.
        \end{itemize}
    \item We argue that the $\suspindex$-th horn of $\widetilde{ x}$ defines a map of simplicial sets
    \[\Lambda^{\suspindex}[m+1]\to W_{\suspindex-1}\]
    with  marking.

    Let $[\mathbf{a''}]$ be a marked $p$-simplex of $\Lambda^\suspindex[m+1]$, namely $[\mathbf{a}]$
    contains the vertices $\{\suspindex-1, \suspindex, \suspindex+1\}\cap [m+1]$. If the simplex $\widetilde x([\mathbf{a''}])$ is totally degenerate, it is in particular marked, so we will exclude this case for the rest of the discussion. 
    By definition of the suspect index, we have $0\leq \suspindex\leq k+1$.  
    Note that $\suspindex=0$ would imply $x=s_{k+1}d_{k+1}x$ using that $\ell >0$, the characterization \cref{CriterionDeg} and \ref{CondAcross}.
    Thus, we can assume $0< \suspindex\leq k+1<k+1+\ell=m+1$. In particular, $\{\suspindex-1, \suspindex, \suspindex+1\}\subseteq [m+1]$.

    If $r\leq k$, using \ref{SuspRowAcross} or \ref{SuspColumnAcross} we have
    \[\widetilde x([\mathbf{a''}])=\widetilde{x}([\mathbf{a'},r-1,r,r+1,\mathbf{a}]\otimes [\mathbf{b}])=0,\]
    so $\widetilde x([\mathbf{a''}])$ is marked by \cref{nerveofsuspension}.
    If instead $r=k+1$, then using \ref{CondAcross} we obtain
    \[\widetilde{x}([\mathbf{a'},r-1,r]\otimes [0,\mathbf{b}])=0,\]
     so $\widetilde x([\mathbf{a''}])$ is marked by \cref{nerveofsuspension}.
    In total, we see that such a face is necessarily marked. Since this holds for all faces as above, we indeed obtain a map of simplicial sets with marking 
     \[\Lambda^{\suspindex}[m+1]\to W_{\suspindex-1}.\]

    \item If furthermore $ x$ is marked, we argue that the $\suspindex$-th horn of $\widetilde{ x}$ defines a map of simplicial sets with  marking
        \[\Lambda^{\suspindex}[m+1]'\to W_{\suspindex-1},\]
        with the simplicial set with marking $\Lambda^{\suspindex}[m+1]'$ defined in \cref{CompMarkAtOnce}.

        We show that the $(\suspindex-1)$-st face is marked using \cref{nerveofsuspension}. 
If $\suspindex\leq k$, since $\ell\geq 1$ we can use \ref{SuspColumnAcross} and obtain
        \[
     \widetilde{x}([0,\ldots, \widehat{\suspindex-1}, \ldots, k+1]\otimes [0,\ldots, \ell])=0,
        \]
        so the $(r-1)$-st face is marked in this case. 
    If $\suspindex=k+1$, by \ref{SuspCritRow} or \ref{SuspArrayAbove} we obtain
        \[
      \widetilde{x}([0,\ldots, \widehat{k}, k+1]\otimes [0,\ldots, \ell])=0
        \]
        so the $(r-1)$-st face is marked in this case.
        
        We show that the $(\suspindex+1)$-st face is marked using \cref{nerveofsuspension}.
        If $\suspindex\leq k$, since $\ell\geq 1$, we can use \ref{SuspArrayAbove} or \ref{SuspRowAcross} to obtain that
        \[
      \widetilde{x}([0,\ldots, \widehat{\suspindex+1}, \ldots, k+1]\otimes [0,\ldots, \ell])=0.
        \]
        so the $(k+1)$-st face is marked in this case.
If $\suspindex=k+1$, using \ref{SuspArrayAbove} and the fact that $x$ is marked we obtain
        \[
      \widetilde{x}([0,\ldots, k+1]\otimes [1,\ldots, \ell])=
        x([0,\ldots, k]\otimes [0,1,\ldots, \ell])=0.
        \]
So the $(k+1)$-st face is marked in this case.

\end{itemize}
 By filling all $\suspindex$-horns of suspect $(m+1)$-simplices $\widetilde{ x}$ of $W_{\suspindex}$, we then obtain their $\suspindex$-th face $ x$, which was missing in $W_{\suspindex-1}$, as well as the suspect $(m+1)$-simplex $\widetilde{ x}$ itself.
This can be rephrased by saying that there is a pushout square
\[
\begin{tikzcd}[column sep=1.0cm]
\underset{\mathclap{\begin{subarray}{c}
   x \\
  \mbox{\tiny{non-marked}}
  \end{subarray}} }{\coprod} \Lambda^{\suspindex}[m+1]\amalg\underset{\mathclap{\begin{subarray}{c}
  x \\
  \mbox{\tiny{marked}}
  \end{subarray} }}{\coprod} \Lambda^{\suspindex}[m+1]'\arrow[d]\arrow[r] &
 \underset{\mathclap{\begin{subarray}{c}
   x \\
  \mbox{\tiny{non-marked}}
  \end{subarray} }}{\coprod} \Delta^{\suspindex}[m+1]\amalg \underset{\mathclap{\begin{subarray}{c}
 x\\
  \mbox{\tiny{marked}}
  \end{subarray}} }{\coprod} \Delta^{\suspindex}[m+1]''\arrow[d]\\
W_{\suspindex-1} \arrow[r]& W_{\suspindex}.
\end{tikzcd}
\]
Since the involved horn inclusions are in fact inner horn inclusions,
the inclusions of simplicial sets with marking $\Lambda^{\suspindex}[m+1]\hookrightarrow\Delta^{\suspindex}[m+1]$ and $\Lambda^{\suspindex}[m+1]'\hookrightarrow\Delta^{\suspindex}[m+1]''$ are complicial inner anodyne extensions by \cref{CompMarkAtOnce}.

It follows that the inclusion  $W_{{\suspindex-1}}\hookrightarrow W_{\suspindex}$ for any $1\le \suspindex\le m-j$, the inclusion $Y_{j-1}\hookrightarrow Y_j$ for any $1\le j\le m$, the inclusion $Y_{j-1}\hookrightarrow Y_j$ for any $1\le j\le m$, the inclusion $X_{m-1}\hookrightarrow X_m$ for any $m\ge1$, and the inclusion $\Sigma\NRS\cC\to \NRS\Sigma\cC$ are complicial inner anodyne extensions, as desired.
\end{proof}

\section{$\Theta_n$-spaces and Quillen pair with complicial sets}

In this section we apply \cref{ThmB} to produce an explicit Quillen adjunction between the model structure of $n$-complicial sets, and the model structure for complete Segal $\Theta_n$-spaces, which we first recall.

\subsection{$\Theta_n$-spaces}

We recall the main facts about $\Theta_n$-spaces that will be used in this paper. 

\begin{rmk}
Let $n\ge0$.
The suspension functor $\Sigma\colon\omega\cat\to\omega\cat_{*,*}$ restricts and corestricts to a functor $\Sigma\colon(n-1)\cat\to n\cat_{*,*}$.
\end{rmk}

Let $\Theta_n$ denote Joyal's cell category from \cite{JoyalDisks}, which is by \cite{BergerCellular,MakkaiZawadowski} a full subcategory of $n\cat$. By definition, $\Theta_0$ is the terminal category, $\Theta_1$ is the ordinal category $\Delta$, and
$\Theta_n$ is for $n>0$ the full subcategory of $n\cat$ whose generic object is obtained as a pushout of $n$-categories
\[\theta=[k|\theta_1,\dots,\theta_k]=\Sigma\theta_1\aamalg{[0]}\Sigma\theta_2\aamalg{[0]}\dots\aamalg{[0]}\Sigma\theta_k\]
for $k\ge0$ and $\theta_1,\dots,\theta_k\in\Theta_{n-1}$.

For $n>0$, there is a full inclusion $\Theta_{n-1}\hookrightarrow\Theta_n$, and whenever needed we will regard any object of $\Theta_{n-1}$ as an object of $\Theta_n$ without further specification.

\begin{defn}
Let $\ge0$.
A \emph{$\Theta_n$-space} (resp.~\emph{$\Theta_n$-set}) is a presheaf $W\colon\Theta_n^{\op}\to\sset$ (resp.~$W\colon\Theta_n^{\op}\to\set$).
\end{defn}

For $n\ge0$, we denote by $\spsh{\Theta_n}$ (resp. $\psh{\Theta_n}$) the category of $\Theta_n$-spaces (resp.~$\Theta_n$-sets).

\begin{rmk}
The categories $\spsh{\Theta_n}$ and $\psh{\Theta_n}$ are cocomplete, and colimits are computed componentwise.
\end{rmk}

For $n\ge0$ the canonical inclusion $\set\hookrightarrow\sset$ of sets as discrete simplicial sets induces a canonical inclusion $\psh{\Theta_n}\hookrightarrow\spsh{\Theta_n}$, which preserves limits and colimits. In particular, we often regard $\Theta_n$-sets as discrete $\Theta_n$-spaces without further specification.

For any object $\theta$ in $\Theta_n$, we denote by $\Theta_n[\theta]$ the $\Theta_n$-set represented by $\theta$ via the Yoneda embedding $\Theta_n\hookrightarrow\psh{\Theta_n}$.

\begin{rmk}\label{BoxProduct}
As a special case of \cite[\textsection 3.1]{Ara}, given any $\Theta_n$-set $A$ and any space $B$ one can consider the $\Theta_n$-space $A\boxtimes B$, which is defined levelwise as the simplicial set
\[(A\boxtimes B)_{\theta}:=A_{\theta}\times B.\]
The construction extends to a bifunctor
\[\boxtimes\colon\psh{\bT}\times\sset\to\spsh{\bT}\]
that preserves colimits in each variable.
\end{rmk}

\subsection{Suspension of $\bT$-spaces}

\begin{rmk}
Let $n\ge0$.
The suspension functor $\Sigma\colon\omega\cat\to\omega\cat$ restricts and corestricts to a functor $\Sigma\colon\Theta_{n-1}\to \Theta_n$.
\end{rmk}

As discussed in \cite[Remark~4.5]{rezkTheta} and in \cite[Lemma~11.10]{rezkTheta}, the following functor agrees with the functor $V[1]$ from \cite[\textsection~4.4]{rezkTheta}.

\begin{defn}
Let $n>0$, and $\theta\in\Theta_{n-1}$. The \emph{suspension} of the representable presheaf $\Theta_{n-1}[\theta]$ is the (discrete) $\Theta_n$-space
\[\Sigma\Theta_n[\theta]\coloneqq\Theta_n[\Sigma\theta].\]
The enriched left Kan extension of the functor $\Theta_{n-1}\to\Theta_n\hookrightarrow\spsh{\Theta_n}$ defines a functor $\Sigma\colon\spsh{\Theta_{n-1}}\subseteq\spsh{\Theta_{n}}_{*,*}$.
\end{defn}

\subsection{The adjunction}

Let us begin by defining the functor that we use to make our comparison.

\begin{const}
Let $n \ge0$.
The functor $\bT\times\Delta\subseteq\spsh{\bT}\to\msset$ given by 
\[(\theta,[\ell])\mapsto(\bT\times\Delta)[\theta,\ell]=\bT[\theta]\boxtimes\Delta[\ell]\mapsto \NRS\theta\times\Delta[\ell]^\sharp.\]
induces an adjunction
\[L_n\colon\spsh{\bT}\rightleftarrows\msset\colon R_n.\]
\end{const}

\subsection{The Quillen pair before localizing}

\begin{prop}
\label{projmodelstructure}
Let $n\ge0$.
The category $\spsh{\Theta_n}$ supports the projective
model structure $\spsh{\Theta_n}_{(\infty,n)}$ where the fibrant objects are precisely the projectively fibrant presheaves and the cofibrations are precisely the projective cofibrations.
\end{prop}

\begin{rmk}
Let $n\ge0$.
Combining \cite[Theorem 11.6.1, Definition. 11.5.33, Definition 11.5.25]{hirschhorn}, we know that
\begin{enumerate}[leftmargin=*]
    \item a set of  generating cofibrations for the projective model structure on $\spsh{\bT}$ is given by all maps of the form
    \[\bT[\theta]\boxtimes\partial\Delta[\ell]\to\bT[\theta]\boxtimes\Delta[\ell]\quad\text{ for }\theta \in \bT \mbox{ and }\ell\geq 0;\]
    and 
    
    \item a set of generating acyclic cofibrations for the projective model structure on $\spsh{\bT}$ is given by all maps of the form
    \[\bT[\theta]\boxtimes\Lambda^k[\ell]\to\bT[\theta]\boxtimes\Delta[\ell]\quad\text{ for }\theta \in \bT \mbox{ and } 0\leq k \leq \ell.\]
  \end{enumerate}
\end{rmk}

The following can be proven similarly to {\cite[Lemma~1.27]{BOR}}.

\begin{prop}
\label{SharpLeftQuillen}
Let $\ge0$.
The functor
\[(-)^{\sharp}\colon\sset_{(\infty,0)}\to \msset_{(\infty,n)}\] is left Quillen.
\end{prop}

\begin{prop}
Let $n\ge0$.\label{Linjective}
The functor
\[L_n\colon\spsh{\bT}_{\mathrm{proj}}\to\msset_{(\infty,n)}\]
is left Quillen.
\end{prop}

We include the proof for the reader's convenience, but the argument is the evident generalization of the $2$-dimensional case treated in \cite[Proposition~2.2]{BOR}.

\begin{proof}
We want to show that the functor $L_n$ preserves cofibrations and acyclic cofibrations.
Using the facts that $(-)^{\sharp}$ commutes with colimits, which is a consequence of \cref{SharpLeftQuillen}, and that the box product $\boxtimes$ preserves colimits in each variable, which was recalled in \cref{BoxProduct}, we see that
\begin{enumerate}[leftmargin=*]
\item the image of the generic generating cofibration via $L_n$ is the map
 \[\NRS\theta\times\partial\Delta[\ell]^\sharp\to\NRS\theta\times\Delta[\ell]^\sharp\quad\text{ for }\theta \in \bT \mbox{ and }\ell\geq 0;\]
 
\item the image of the generic generating acyclic cofibration via $L_n$ is the map
\[\NRS\theta\times\Lambda^k[\ell]^\sharp\to\NRS\theta\times\Delta[\ell]^\sharp\quad\text{ for }\theta \in \bT \mbox{ and } 0\leq k \leq \ell.\]
    \end{enumerate}
    Since the model structure $\msset_{(\infty,n)}$ is cartesian closed by \cref{modelstructurewithsaturation} and $(-)^{\sharp}$ is a left Quillen functor by \cref{SharpLeftQuillen}, we conclude that
    \begin{enumerate}[leftmargin=*]
    \item the map $\NRS\theta\times\partial\Delta[\ell]^\sharp\to\NRS\theta\times\Delta[\ell]^\sharp$ is a cofibration and
    \item the map $\NRS\theta\times\Lambda^k[\ell]^\sharp\to\NRS\theta\times\Delta[\ell]^\sharp$ is an acyclic cofibration
    \end{enumerate}
    It follows that $L_n$ preserves cofibrations and acyclic cofibrations, so it is a left Quillen functor, as desired.
\qedhere
\end{proof}

\subsection{The Quillen pair before localizing}

We construct a variant of Rezk's model structure from \cite[\textsection11.4]{rezkTheta}  (see \cref{changesfromRezk}).

\begin{defn}
\label{anodynemapstheta}
Let $n\ge0$.
A map of (discrete) $\Theta_n$-spaces is an \emph{elementary acyclic cofibration} if it is of one of the following kinds:
\begin{enumerate}[leftmargin=*]
    \item For $2\leq j \leq n$, $k\ge 1$ and objects $\theta_1,\ldots, \theta_k$ of $\Theta_{n-j}$, the \emph{$j$-fold $k$-Segality extension}
        \[ \Sigma^j\bT[\theta_1]\aamalg{\Sigma^{j-1}\bT[0]}\dots\aamalg{\Sigma^{j-1}\bT[0]}\Sigma^j\bT[\theta_k]\hookrightarrow\Sigma^{j-1}\bT[k|\theta_1, \ldots, \theta_k]\]
        \item For $0\leq j \leq n-1$, the \emph{$j$-fold completeness extension}
    \[ \Sigma^j\bT[0]\hookrightarrow\Sigma^j\bT[0]\aamalg{\Sigma^j\bT[1]}\Sigma^j\bT[3]\aamalg{\Sigma^j\bT[1]}\Sigma^j\bT[0].\]
\end{enumerate}
\end{defn}

\begin{defn}
\label{defcomplicial}
A \emph{complete Segal $\Theta_n$-space}
is a $\Theta_n$-space that is local with respect to all maps of type (1) and (2) from \cref{anodynemapstheta}.
\end{defn}

By localizing the projective model structure $\psh{\Theta_n}_{\mathrm{proj}}$ at the class of maps from \cref{anodynemapstheta}, we obtain the following.

\begin{thm}
\label{modelstructurelocalized}
Let $n\ge0$.
\label{modelstructureondiscretepresheaves}
The category $\spsh{\Theta_n}$ supports a cartesian closed
model structure $\spsh{\Theta_n}_{(\infty,n)}$ where the fibrant objects are precisely the projectively fibrant complete Segal $\Theta_n$-spaces and the cofibrations are precisely the projective cofibrations.
\end{thm}

\begin{rmk}
\label{changesfromRezk}
The model structure from \cref{modelstructurelocalized} differs from Rezk's from \cite[\textsection11.4]{rezkTheta} in the following aspects:
\begin{enumerate}[leftmargin=*]
\item We work with localizations of the projective model structure, instead of the injective model structure.
\item To express the completeness extension, we use $\bT[0]\amalg_{\bT[1]}\bT[3]\amalg_{\bT[1]}\bT[0]$ instead of the $\Theta_n$-nerve of the walking isomorphism.
\end{enumerate}
However, the two model structures are Quillen equivalent (see \cite[\textsection2.5--2.13,\textsection10]{rezkTheta}).
\end{rmk}

We now show that we still have a Quillen pair after localizing the projective model structure on $\spsh{\bT}$.

\begin{thm} \label{LQuillenPair}
Let $n\ge0$.
The functor
\[L_n\colon\spsh{\bT}_{(\infty,n)}\to\msset_{(\infty,n)}\]
is left Quillen.
\end{thm}

\begin{proof}
We prove this by induction on $n\ge2$. The basis of the induction $n=0,1,2$ are treated in \cref{SharpLeftQuillen}, the combination of \cite[\textsection6.5]{VerityComplicialI} with \cite{JT}, and \cite[Theorem~2.4]{BOR}, respectively. We now assume $n>2$ and that the statement is true for $n-1$.
Since cofibrations are unchanged by localization, it suffices to prove that $L_n$ preserves acyclic cofibrations.  We do so by proving that $L_n$ preserves all elementary acyclic cofibrations from \cref{anodynemapstheta}, which we do in \cref{Segality,Completeness}.
\end{proof}

The following result by Steiner will allow us to apply \cref{ThmB} to the case $\cC=\theta$, for some $\theta\in\Theta_n$.

\begin{thm}[{\cite{SteinerSimpleOmega}}]
\label{thetaalgebraic}
Let $n\ge0$ and $\theta\in\Theta_n$. There is an isomorphism of $\omega$-categories
\[\theta\cong\nu T\]
for some augmented directed chain complex $T$.
\end{thm}

We analyze the action of $L_n$ on the Segality extensions.

\begin{prop}
\label{Segality}
Let $n>0$.
If the functor
$L_{n-1}\colon\spsh{\bT}_{(\infty,n-1)}\to\msset_{(\infty,n-1)}$
is left Quillen, the functor $L_n\colon\spsh{\Theta_{n}}\to\msset$ sends the $j$-fold Segal acyclic cofibration from \cref{anodynemapstheta},
\[ \Sigma^j\bT[\theta_1]\aamalg{\Sigma^{j-1}\bT[0]}\dots\aamalg{\Sigma^{j-1}\bT[0]}\Sigma^j\bT[\theta_k]\hookrightarrow\Sigma^{j-1}\bT[k|\theta_1,\dots,\theta_k] \]
for $1\leq j \leq n$, $k\ge 1$ and objects $\theta_1,\ldots, \theta_k$ of $\Theta_{n-j}$,
  to a weak equivalence in $\msset_{(\infty,n)}$.
\end{prop}

\begin{proof}
We prove the statement by induction on $j\ge1$, and fixed $k\ge 1$. The basis case(s) $j=1$ is a direct consequence of \cite[Theorem~4.9]{ORfundamentalpushouts}, and we now assume $j>1$ for the inductive step.
By definition, the map
  \[\Sigma^{j-1}\bT[\theta_1]\aamalg{\Sigma^{j-2}\bT[0]}\dots\aamalg{\Sigma^{j-2}\bT[0]}\Sigma^{j-1}\bT[\theta_k]\hookrightarrow\Sigma^{j-2}\bT[k|\theta_1,\dots,\theta_k].\]
  is an acyclic cofibration in $\spsh{\Theta_{n-1}}_{(\infty,n-1)}$.
  This acyclic cofibration can be rewritten as
  \[ \bT\Sigma^{j-1}[\theta_1]\aamalg{\bT\Sigma^{j-2}[0]}\dots\aamalg{\bT\Sigma^{j-2}[0]}\bT\Sigma^{j-1}[\theta_k]\hookrightarrow\bT\Sigma^{j-2}[k|\theta_1,\dots,\theta_k].\]
By applying to it the left Quillen functor $L_{n-1}\colon\spsh{\bT}_{(\infty,n-1)}\to\msset_{(\infty,n-1)}$ we obtain an acyclic cofibration in $\msset_{(\infty,n-1)}$
  \[\NRS\bT\Sigma^{j-1}[\theta_1]\aamalg{ \NRS\bT\Sigma^{j-2}[0]}\dots\aamalg{ \NRS\bT\Sigma^{j-2}[0]} \NRS\bT\Sigma^{j-1}[\theta_k]\hookrightarrow \NRS\bT\Sigma^{j-2}[k|\theta_1,\dots,\theta_k].\]
By applying to it the left Quillen functor $\Sigma\colon\msset_{(\infty,n-1)}\to(\msset_{(\infty,n)})_{*,*}$ from \cref{suspensionhomotopical} we obtain an acyclic cofibration in $\msset_{(\infty,n)}$
   \[ \Sigma\NRS\bT\Sigma^{j-1}[\theta_1]\aamalg{ \Sigma\NRS\bT\Sigma^{j-2}[0]}\dots\aamalg{ \Sigma\NRS\bT\Sigma^{j-2}[0]} \Sigma\NRS\bT\Sigma^{j-1}[\theta_k]\hookrightarrow \Sigma\NRS\bT\Sigma^{j-2}[k|\theta_1\dots,\theta_k].\]
Since $\Sigma$ commutes with nerve by \cref{ThmB,thetaalgebraic}, we also get an acylic cofibration
  \[ \NRS\bT\Sigma^{j}[\theta_1]\aamalg{ \NRS\bT\Sigma^{j-1}[0]}\dots\aamalg{ \NRS\bT\Sigma^{j-1}[0]} \NRS\bT\Sigma^{j}[\theta_k]\hookrightarrow \NRS\bT\Sigma^{j-1}[k|\theta_1\dots,\theta_k],\]
   which is
    \[ L_n\bT\Sigma^{j}[\theta_1]\aamalg{ L_n\bT\Sigma^{j-1}[0]}\dots\aamalg{ L_n\bT\Sigma^{j-1}[0]} L_n\bT\Sigma^{j}[\theta_k]\hookrightarrow L_n\bT\Sigma^{j-1}[k|\theta_1\dots,\theta_k].\]
This concludes the proof.
  \end{proof}

We analyze the action of $L_n$ on the completeness extensions.

\begin{prop} \label{Completeness}
Let $n>0$.
If the functor
$L_{n-1}\colon\spsh{\bT}_{(\infty,n-1)}\to\msset_{(\infty,n-1)}$
is left Quillen, the functor $L_n\colon\spsh{\Theta_{n}}\to\msset$ sends the $j$-fold Segal acyclic cofibration from \cref{anodynemapstheta},
 \[ \Sigma^j\bT[0]\hookrightarrow\Sigma^j\bT[0]\aamalg{\Sigma^j\bT[1]}\Sigma^j\bT[3]\aamalg{\Sigma^j\bT[1]}\Sigma^j\bT[0]\]
 for $0\leq j \leq n-1$,
  to a weak equivalence in $\msset_{(\infty,n)}$.
\end{prop}

\begin{proof}
We prove the statement by induction on $j\ge0$. The basis case(s) $j=0,1$ are proven in \cite[Propositions~2.7,~2.9]{BOR}, and we now assume $j>0$ for the inductive step.
By definition, the map
  \[ \Sigma^{j-1}\bT[0]\hookrightarrow\Sigma^{j-1}\bT[0]\aamalg{\Sigma^{j-1}\bT[1]}\Sigma^{j-1}\bT[3]\aamalg{\Sigma^{j-1}\bT[1]}\Sigma^{j-1}\bT[0]\]
  is an acyclic cofibration in $\spsh{\Theta_{n-1}}_{(\infty,n-1)}$. This acyclic cofibration can be rewritten as
   \[ \bT\Sigma^{j-1}[0]\hookrightarrow\bT\Sigma^{j-1}[0]\aamalg{\bT\Sigma^{j-1}[1]}\bT\Sigma^{j-1}[3]\aamalg{\bT\Sigma^{j-1}[1]}\bT\Sigma^{j-1}[0].\]
By applying to it the left Quillen functor
$L_{n-1}\colon\spsh{\bT}_{(\infty,n-1)}\to\msset_{(\infty,n-1)}$
we obtain an acyclic cofibration in $\msset_{(\infty,n-1)}$
   \[ \NRS\bT\Sigma^{j-1}[0]\hookrightarrow\NRS\bT\Sigma^{j-1}[0]\aamalg{\NRS\bT\Sigma^{j-1}[1]}\NRS\bT\Sigma^{j-1}[3]\aamalg{\NRS\bT\Sigma^{j-1}[1]}\NRS\bT\Sigma^{j-1}[0].\]
   By applying to it the left Quillen functor
$\Sigma\colon\msset_{(\infty,n-1)}\to(\msset_{(\infty,n)})_{*,*}$ from \cref{suspensionhomotopical} we obtain an acyclic cofibration in $(\msset_{(\infty,n)})_{(*,*)}$
   \[ \Sigma\NRS\bT\Sigma^{j-1}[0]\hookrightarrow\Sigma\NRS\bT\Sigma^{j-1}[0]\aamalg{\Sigma\NRS\bT\Sigma^{j-1}[1]}\Sigma\NRS\Sigma^{j-1}\bT[3]\aamalg{\Sigma\NRS\bT\Sigma^{j-1}[1]}\Sigma\NRS\bT\Sigma^{j-1}[0].\]
Since $\Sigma$ commutes with nerve by \cref{ThmB,thetaalgebraic} we also get an acyclic cofibration
   \[ \NRS\bT\Sigma^{j}[0]\hookrightarrow\NRS\bT\Sigma^{j}[0]\aamalg{\NRS\bT\Sigma^{j}[1]}\NRS\bT\Sigma^{j}[3]\aamalg{\NRS\bT\Sigma^{j}[1]}\NRS\bT\Sigma^{j}[0],\]
   which is
      \[L_n\bT\Sigma^{j}[0]\hookrightarrow L_n\bT\Sigma^{j}[0]\aamalg{L\bT\Sigma^{j}[1]}L_n\bT\Sigma^{j}[3]\aamalg{L_n\bT\Sigma^{j}[1]}L_n\bT\Sigma^{j}[0].\]
This concludes the proof.
  \end{proof}

With the establishment of \cref{Segality,Completeness}, the proof of \cref{LQuillenPair} is now complete.

\appendix

\label{appendix}

\section{Proof of \cref{parent}}

We now prove \cref{parent}.

\begin{proof}[Proof of \cref{parent}]
Since these cases are mutually exclusive and cover all possibilities, this at least defines a map, and by construction we will have $d_{\suspindex}\widetilde x=x$. 

It is immediate that the map is directed. Observe that neither the case \ref{SuspColumnAbove}
nor the 
case \ref{SuspColumnAcross} can apply to a chain of dimension $0$, proving that $\widetilde x$ is augmented since $x$ is augmented. What we need to check is that $\widetilde x$ is a chain map.

\begin{enumerate}[leftmargin=*, label= \normalfont{(P\arabic*)}]
    \item This case is immediate since $x$ is a chain map (and \lq{}not containing $\suspindex$ in the first
    component\rq{} is preserved by the differential).
    \item For $\lvert \mathbf{a} \rvert=-1$, there is nothing to check. For $\lvert \mathbf{a} \rvert\geq 0$, on the one hand, we have 
    \[
    \begin{array}{llll}
     \widetilde x(\partial [\suspindex,\mathbf{a}]\otimes [b])    &=& \widetilde x([\mathbf{a}]\otimes[b]-[\suspindex,\partial\mathbf{a}]\otimes [b] )&\\
    &=& x(s^{\suspindex-1}[\mathbf{a}]\otimes[b])-x(s^{\suspindex-1}[\suspindex,\partial\mathbf{a}]\otimes[0]);&
     \end{array}
    \]
    on the other hand we have
        \[
    \begin{array}{llll}
     \partial \widetilde x([\suspindex,\mathbf{a}]\otimes [b])&=&\partial x(s^{\suspindex-1}[\suspindex,\mathbf{a}]\otimes[0])&\\
     &=&x(s^{\suspindex-1}[\mathbf{a}]\otimes[0])-x(s^{\suspindex-1}[\suspindex,\partial\mathbf{a}]\otimes[0]).&
    \end{array}
    \]
From \ref{CondBelow}, we obtain
    \[
    \begin{array}{llll}
    0=\partial x(s^{\suspindex-1}[\mathbf{a}]\otimes[0,b]) &=& x(s^{\suspindex-1}[\partial\mathbf{a}]\otimes[0,b]) &\\
    &&+(-1)^{\lvert \mathbf{a} \rvert}x(s^{\suspindex-1}[\mathbf{a}]\otimes[0])\\
    &&-(-1)^{\lvert \mathbf{a} \rvert}x(s^{\suspindex-1}[\mathbf{a}]\otimes[b]).&
    \end{array}
    \]
  Using \ref{CondBelow} again, the first summand vanishes, yielding the equality of the other two. This shows the desired equality.

    \item For $\lvert \mathbf{a'}\rvert\geq 1$, $\lvert \mathbf{b}\rvert \geq 2$,  on the one hand we have
    \[
    \begin{array}{llll}
     \widetilde x(\partial ([\mathbf{a'},\suspindex] \otimes [\mathbf{b}]))    &=& \widetilde x([\partial \mathbf{a'},\suspindex] \otimes [\mathbf{b}]\\
     &&     + (-1)^{\lvert \mathbf{a'}\rvert+1}[\mathbf{a'}] \otimes [\mathbf{b}]+ (-1)^{\lvert \mathbf{a'}\rvert+1}[\mathbf{a'},\suspindex] \otimes \partial^{\circ}[\mathbf{b}])& \\
         &=&x([\partial\mathbf{a'}] \otimes [0,\mathbf{b}])+ (-1)^{\lvert \mathbf{a'}\rvert+1}x([\mathbf{a'}] \otimes [\mathbf{b}])+&\\
         &&+(-1)^{\lvert \mathbf{a'}\rvert+1}x([\mathbf{a'}] \otimes [0,\partial^{\circ}\mathbf{b}]);&
    \end{array}
    \]
    on the other hand, we have 
      \[
    \begin{array}{llll}
     \partial \widetilde x( [\mathbf{a'},\suspindex] \otimes [\mathbf{b}])    &=& \partial(x([\mathbf{a'}] \otimes [0,\mathbf{b}]))& \\
          &=&x([\partial\mathbf{a'}] \otimes [0,\mathbf{b}])\\
          &&+ (-1)^{\lvert \mathbf{a'}\rvert+1}x([\mathbf{a'}] \otimes [\mathbf{b}])&\\
         &&+(-1)^{\lvert \mathbf{a'}\rvert+1}x([\mathbf{a'}] \otimes [0,\partial^{\circ}\mathbf{b}]),&
    \end{array}
    \]
so the two expressions coincide.

   For $\lvert\mathbf{a'}\rvert=0$, $\lvert \mathbf{b}\rvert \geq 2$, on the one hand we have
    \[
     \begin{array}{llll}
     \widetilde x(\partial ([a',\suspindex] \otimes [\mathbf{b}]))    &=& \widetilde x([\suspindex] \otimes [\mathbf{b}]-[a'] \otimes [\mathbf{b}]-[a',\suspindex] \otimes \partial^{\circ}[\mathbf{b}])& \\
         &=&-x([a'] \otimes [\mathbf{b}])-x([a'] \otimes [0,\partial^{\circ}\mathbf{b}]);&

    \end{array}
    \]
    On the other hand we have 
      \[
    \begin{array}{llll}
     \partial \widetilde x( [a',\suspindex] \otimes [\mathbf{b}])    &=& \partial(x([a'] \otimes [0,\mathbf{b}]))& \\
          &=& -x([a'] \otimes [\mathbf{b}])-x([a'] \otimes [0,\partial^{\circ}\mathbf{b}]),&
    \end{array}
    \]
so the two expressions coincide.

   For $\lvert \mathbf{a'}\rvert\geq 1$, $\lvert \mathbf{b}\rvert =1$, on the one hand we have
    \[
    \begin{array}{llll}
     \widetilde x(\partial ([\mathbf{a'},\suspindex] \otimes [b_0,b_1]))    &=& \widetilde x([\partial \mathbf{a'},\suspindex] \otimes [b_0,b_1]\\
     &&+ (-1)^{\lvert \mathbf{a'}\rvert+1}[\mathbf{a'}] \otimes [b_0,b_1]&\\
         &&+ (-1)^{\lvert \mathbf{a'}\rvert+1}[\mathbf{a'},\suspindex] \otimes [b_0]\\
         &&- (-1)^{\lvert \mathbf{a'}\rvert+1}[\mathbf{a'},\suspindex] \otimes [b_1])& \\
         &=&x([\partial\mathbf{a'}] \otimes [0,b_0,b_1])\\
         &&+ (-1)^{\lvert \mathbf{a'}\rvert+1}x([\mathbf{a'}] \otimes [b_0,b_1])&\\
         &&+(-1)^{\lvert \mathbf{a'}\rvert+1}x([\mathbf{a'}] \otimes [0,b_0])\\
         &&+x([\mathbf{a'},\suspindex-1] \otimes [0])&\\
         &&-(-1)^{\lvert \mathbf{a'}\rvert+1}x([\mathbf{a'}] \otimes [0,b_1])\\
         &&-x([\mathbf{a'},\suspindex-1] \otimes [0]);
        &
    \end{array}
    \]
    on the other hand we have 
      \[
    \begin{array}{llll}
     \partial \widetilde x( [\mathbf{a'},\suspindex] \otimes [b_0,b_1])    &=& \partial(x([\mathbf{a'}] \otimes [0,b_0,b_1]))& \\
          &=&x([\partial\mathbf{a'}] \otimes [0,b_0,b_1])\\
          &&+ (-1)^{\lvert \mathbf{a'}\rvert+1}x([\mathbf{a'}] \otimes [b_0,b_1])&\\
         &&+(-1)^{\lvert \mathbf{a'}\rvert}x([\mathbf{a'}] \otimes [0,b_1])\\
         &&-(-1)^{\lvert \mathbf{a'}\rvert}x([\mathbf{a'}] \otimes [0,b_0]),&
    \end{array}
    \]
    so the two expressions coincide.
    
For $\lvert \mathbf{a'}\rvert=0$, $\lvert \mathbf{b}\rvert =1$, on the one hand we have
    \[
    \begin{array}{llll}
     \widetilde x(\partial ([a',\suspindex] \otimes [b_0,b_1]))    &=& \widetilde x([r] \otimes [b_0,b_1]-[a'] \otimes [b_0,b_1]&\\
         &&- [a',\suspindex] \otimes [b_0]+[a',\suspindex] \otimes [b_1])& \\
         &=&-x[a'] \otimes [b_0,b_1]-x([a'] \otimes [0,b_0])&\\
         &&-x([a',\suspindex-1] \otimes [0])+x([a'] \otimes [0,b_1])&\\
         &&+x([a',\suspindex-1] \otimes [0]);&
    \end{array}
    \]
    on the other hand we have 
      \[
    \begin{array}{llll}
     \partial \widetilde x( [\mathbf{a'},\suspindex] \otimes [b_0,b_1])    &=& \partial(x([a'] \otimes [0,b_0,b_1]))& \\
          &=&-x([a'] \otimes [b_0,b_1]-x([a'] \otimes [0,b_0])\\
          &&+x([a'] \otimes [0,b_1]),
    \end{array}
    \]
so the two expressions coincide.

    \item For $\lvert \mathbf{a'}\rvert\geq 1$, on the one hand we have
      \[
    \begin{array}{llll}
     \widetilde x(\partial ([\mathbf{a'},\suspindex] \otimes [b]) )   &=& \widetilde x((-1)^{\lvert \mathbf{a'}\rvert+1}[\mathbf{a'}]\otimes [b] + [(\partial \mathbf{a'}),\suspindex]\otimes [b])& \\
         & =&(-1)^{\lvert \mathbf{a'}\rvert+1}x([\mathbf{a'}]\otimes [b])+x([(\partial \mathbf{a'})]\otimes [0,b])\\
         &&+x([(\partial \mathbf{a'}),\suspindex-1]\otimes [0]);
    \end{array}
    \]
      on the other hand we have 
      \[ 
    \begin{array}{llll}
     \partial \widetilde x( [\mathbf{a'},\suspindex] \otimes [b])    &=& \partial (x([\mathbf{a'})]\otimes [0,b])+x ([\mathbf{a'},\suspindex-1] \otimes [0]))& \\
        & =&x([(\partial \mathbf{a'})]\otimes [0,b])-(-1)^{\lvert \mathbf{a'}\rvert}x([\mathbf{a'}]\otimes [b])\\
        &&+(-1)^{\lvert \mathbf{a'}\rvert}x([\mathbf{a'}]\otimes [0])+x([(\partial \mathbf{a'}),\suspindex-1]\otimes [0])\\
        &&+(-1)^{\lvert \mathbf{a'}\rvert+1}x([\mathbf{a'}]\otimes [0]);&
    \end{array}
    \]
so the two expressions coincide.

For $\lvert \mathbf{a'}\rvert = 0$,
    on the one hand we have
      \[
    \begin{array}{llll}
     \widetilde x(\partial ([a',\suspindex] \otimes [b]) )   &=& \widetilde x([\suspindex]\otimes [b] - [a']\otimes [b])& \\
         & =&x([\suspindex-1]\otimes [0])-x([a']\otimes [b]);
    \end{array}
    \]
          on the other hand we have 
      \[
    \begin{array}{llll}
     \partial \widetilde x( [a',\suspindex] \otimes [b])    &=& \partial (x([a']\otimes [0,b])+F ([a',\suspindex-1] \otimes [0]))& \\
        & =&x([a']\otimes [0]) -x([a']\otimes [b]) \\
        &&+x([\suspindex-1]\otimes[0])- x([a']\otimes[0]);
    \end{array}
    \]
   so the two expressions coincide.
    \item For $\lvert \mathbf{a}\rvert\geq 1, \lvert \mathbf{a'}\rvert\geq 1$,
    on the one hand we have
      \[
    \begin{array}{llll}
     \widetilde x(\partial ([\mathbf{a'},\suspindex,\mathbf{a}] \otimes [b]) )   &=& \widetilde x([\partial\mathbf{a'},\suspindex,\mathbf{a}] \otimes [b])\\
     &&+(-1)^{\lvert \mathbf{a'}\vert+1}\widetilde x([\mathbf{a'},\mathbf{a}]\otimes [b])\\
     &&+(-1)^{\lvert \mathbf{a'}\vert+2}\widetilde x([\mathbf{a'},\suspindex,\partial\mathbf{a}] \otimes [b])&\\
         & =&x(s^{\suspindex-1}[\partial\mathbf{a'},\suspindex,\mathbf{a}] \otimes [0])\\
         &&+(-1)^{\lvert \mathbf{a'}\vert+1}x(s^{\suspindex-1}[\mathbf{a'},\mathbf{a}] \otimes [b])&\\
         & &+(-1)^{\lvert \mathbf{a'}\vert+2}x(s^{\suspindex-1}[\mathbf{a'},\suspindex,\partial\mathbf{a}] \otimes [0]);&
    \end{array}
    \]
      on the other hand we have
      \[
     \begin{array}{llll}
     \partial \widetilde x([\mathbf{a'},\suspindex,\mathbf{a}] \otimes [b])    &=& \partial x(s^{\suspindex-1}[\mathbf{a'},\suspindex,\mathbf{a}] \otimes [0])& \\
          & =&x(s^{\suspindex-1}[\partial\mathbf{a'},\suspindex,\mathbf{a}] \otimes [0])\\
          &&+(-1)^{\lvert \mathbf{a'}\rvert+1}x(s^{\suspindex-1}[\mathbf{a'},\mathbf{a}] \otimes [0])&\\
         & &+(-1)^{\lvert \mathbf{a'}\vert+2}x(s^{\suspindex-1}[\mathbf{a'},\suspindex,\partial\mathbf{a}] \otimes [0]).&
     \end{array}
     \]
      Observe that from \ref{CondAcross}, we obtain
    \[
    \begin{array}{llll}
    0&=&\partial x(s^{\suspindex-1}[\mathbf{a'},\mathbf{a}]\otimes[0,b]) \\
    &=& x(s^{\suspindex-1}[\partial\mathbf{a'},\mathbf{a}]\otimes[0,b]) \\
    &&+ (-1)^{\lvert \mathbf{a'}\rvert+1}x(s^{\suspindex-1}[\mathbf{a'},\partial\mathbf{a}]\otimes[0,b])&\\
    &&+(-1)^{\lvert \mathbf{a'} \rvert+\lvert \mathbf{a}\rvert}x(s^{\suspindex-1}[\mathbf{a'},\mathbf{a}]\otimes[0])\\
    &&-(-1)^{\lvert \mathbf{a'} \rvert+\lvert \mathbf{a}\rvert}x(s^{\suspindex-1}[\mathbf{a'},\mathbf{a}]\otimes[b]).&
    \end{array}
    \]
     The first two summands vanish by \ref{CondAcross} again, yielding the equality of the other two summands. This shows the desired equality.

For $\lvert \mathbf{a}\rvert=0, \lvert \mathbf{a'}\rvert\geq 1$, on the one hand we have
      \[
    \begin{array}{llll}
     \widetilde x(\partial ([\mathbf{a'},\suspindex,a] \otimes [b]) )   &=& \widetilde x([\partial\mathbf{a'},\suspindex,a] \otimes [b])\\
     &&+(-1)^{\lvert \mathbf{a'}\vert+1}\widetilde x([\mathbf{a'},a]\otimes [b])\\
     &&+(-1)^{\lvert \mathbf{a'}\vert+2}\widetilde x([\mathbf{a'},\suspindex] \otimes [b])&\\
         & =&x(s^{\suspindex-1}[\partial\mathbf{a'},\suspindex,a] \otimes [0])\\
         &&+(-1)^{\lvert \mathbf{a'}\vert+1}x(s^{\suspindex-1}[\mathbf{a'},a] \otimes [b])&\\
         & &+(-1)^{\lvert \mathbf{a'}\vert+2}x([\mathbf{a'}] \otimes [0,b])\\
         &&+(-1)^{\lvert \mathbf{a'}\vert+2}x([\mathbf{a'},\suspindex-1] \otimes [0]);
    \end{array}
    \]
      on the other hand we have 
      \[
     \begin{array}{llll}
     \partial \widetilde x([\mathbf{a'},\suspindex,a] \otimes [b])    &=& \partial x(s^{\suspindex-1}[\mathbf{a'},\suspindex,a] \otimes [0])& \\
          & =&x(s^{\suspindex-1}[\partial\mathbf{a'},\suspindex,a] \otimes [0])\\
          &&+(-1)^{\lvert \mathbf{a'}\vert+1}x(s^{\suspindex-1}[\mathbf{a'},a] \otimes [0])&\\
         & &+(-1)^{\lvert \mathbf{a'}\vert+2}x(s^{\suspindex-1}[\mathbf{a'},\suspindex] \otimes [0]).
     \end{array}
     \]
      Observe that from \ref{CondAcross}, we obtain
    \[
    \begin{array}{llll}
    0&=&\partial x(s^{\suspindex-1}[\mathbf{a'},a]\otimes[0,b])\\
    &=& x(s^{\suspindex-1}[\partial\mathbf{a'},a]\otimes[0,b])\\
    &&+ (-1)^{\lvert \mathbf{a'}\rvert+1}x(s^{\suspindex-1}[\mathbf{a'}]\otimes[0,b])&\\
    &&+(-1)^{\lvert \mathbf{a'} \rvert+1}x(s^{\suspindex-1}[\mathbf{a'},a]\otimes[0])\\
    &&-(-1)^{\lvert \mathbf{a'} \rvert+1}x(s^{\suspindex-1}[\mathbf{a'},a]\otimes[b]).&
    \end{array}
    \]
     The first summand vanishes again by \ref{CondAcross}. This implies the desired equality.
     
For $\lvert \mathbf{a}\rvert\geq 1, \lvert \mathbf{a'}\rvert=0$,
     on the one hand we have
      \[
    \begin{array}{llll}
     \widetilde x(\partial ([a',\suspindex,\mathbf{a}] \otimes [b]) )   &=& \widetilde x([\suspindex,\mathbf{a}] \otimes [b])-\widetilde x([a',\mathbf{a}]\otimes [b])\\
     &&+\widetilde x([a',\suspindex,\partial\mathbf{a}] \otimes [b])&\\
         & =&x(s^{\suspindex-1}[\suspindex,\mathbf{a}] \otimes [0])\\
         &&-x(s^{\suspindex-1}[a',\mathbf{a}] \otimes [b])&\\
         & &+x(s^{\suspindex-1}[a',\suspindex,\partial\mathbf{a}] \otimes [0]);&
    \end{array}
    \]
      on the other hand we have 
      \[
     \begin{array}{llll}
     \partial \widetilde x([a',\suspindex,\mathbf{a}] \otimes [b])    &=& \partial x(s^{\suspindex-1}[a',\suspindex,\mathbf{a}] \otimes [0])& \\
          & =&x(s^{\suspindex-1}[\suspindex,\mathbf{a}] \otimes [0])-x(s^{\suspindex-1}[a',\mathbf{a}] \otimes [0])&\\
         & &+x(s^{\suspindex-1}[a',\suspindex,\partial\mathbf{a}] \otimes [0]).
         \\
     \end{array}
     \]
      Observe that from \ref{CondAcross}, we obtain
    \[
    \begin{array}{llll}
    0&=&\partial x(s^{\suspindex-1}[a',\mathbf{a}]\otimes[0,b]) \\
    &=& x(s^{\suspindex-1}[\mathbf{a}]\otimes[0,b]) \\
    &&-x(s^{\suspindex-1}[a',\partial\mathbf{a}]\otimes[0,b])&\\
    &&+(-1)^{\lvert \mathbf{a} \rvert+1}x(s^{\suspindex-1}[a',\mathbf{a}]\otimes[0])\\
    &&-(-1)^{\lvert \mathbf{a} \rvert+1}x(s^{\suspindex-1}[a',\mathbf{a}]\otimes[b]).&
    \end{array}
    \]
     The first summand vanishes again by \ref{CondBelow} and the second by \ref{CondAcross}, so that the other two summands are equal. This implies the desired equality.

For $\lvert \mathbf{a}\rvert=\lvert \mathbf{a'}\rvert=0$, on the one hand we have
      \[
    \begin{array}{llll}
     \widetilde x(\partial ([a',\suspindex,a] \otimes [b]) )   &=& \widetilde x([\suspindex,a] \otimes [b])\\
     &&-\widetilde x([a',a]\otimes [b])\\
     &&+\widetilde x([a',\suspindex] \otimes [b])&\\
         & =&x([s^{\suspindex-1}(\suspindex,a)] \otimes [0])\\
         &&-x(s^{\suspindex-1}[a',a] \otimes [b])&\\
         & &+x(s^{\suspindex-1}[a',\suspindex] \otimes [0])\\
         &&+x([a']\otimes [0,b]);&
    \end{array}
    \]
      on the other hand  we have 
      \[
     \begin{array}{llll}
     \partial \widetilde x([a',\suspindex,a] \otimes [b])    &=& \partial x(s^{\suspindex-1}[a',\suspindex,a] \otimes [0])& \\
          & =&x(s^{\suspindex-1}[\suspindex,a] \otimes [0])\\
          &&-x(s^{\suspindex-1}[a',a] \otimes [0])&\\
         & &+x(s^{\suspindex-1}[a',\suspindex] \otimes [0]).
         \\
     \end{array}
     \]
        Using \ref{CondAcross}, we obtain
    \[
    \begin{array}{llll}
     0&=&\partial x([a',a-1]\otimes [0,b])  \\  &=&x([a-1]\otimes [0,b])\\
     &&-x([a']\otimes [0,b])&  \\
         &&+x([a',a-1]\otimes [b])-x([a',a-1]\otimes [0]).&
    \end{array}
    \]
     Now the first summand vanishes by \ref{CondBelow}. This yields the desired equality.
     
      \item For $\lvert \mathbf{a}\rvert\geq 1, \lvert \mathbf{a'}\rvert\geq 0$ and $\lvert \mathbf{b}\rvert\geq 2$,
    on the one hand we have
      \[
    \begin{array}{llll}
     \widetilde x(\partial ([\mathbf{a'},\suspindex,\mathbf{a}] \otimes [\mathbf{b}]) )   &=& \widetilde x([\partial\mathbf{a'},\suspindex,\mathbf{a}] \otimes [\mathbf{b}])\\
     &&+(-1)^{\lvert \mathbf{a'}\vert+1}\widetilde x([\mathbf{a'},\mathbf{a}]\otimes [\mathbf{b}])\\
     &&+(-1)^{\lvert \mathbf{a'}\vert+2}\widetilde x([\mathbf{a'},\suspindex,\partial\mathbf{a}] \otimes [\mathbf{b}])\\
     &&+(-1)^{\lvert \mathbf{a'}\rvert+2+\lvert \mathbf{a}\rvert}\widetilde x([\mathbf{a'},\suspindex,\mathbf{a}] \otimes \partial^{\circ}[\mathbf{b}])&\\
         & =&(-1)^{\lvert \mathbf{a'}\vert+1}x(s^{\suspindex-1}[\mathbf{a'},\mathbf{a}] \otimes [\mathbf{b}]);&\\
    \end{array}
    \]
      on the other hand we have
      \[
     \begin{array}{llll}
     \partial \widetilde x([\mathbf{a'},\suspindex,\mathbf{a}] \otimes [\mathbf{b}])    &=& 0.&
     \end{array}
     \]
    Observe that from \ref{CondAcross}, we obtain
        \[
    \begin{array}{lllllll}
    0&=&\partial x(s^{\suspindex-1}[\mathbf{a'},\mathbf{a}]\otimes[0,\mathbf{b}])\\
    &=&x(s^{\suspindex-1}[\partial\mathbf{a'},\mathbf{a}]\otimes[0,\mathbf{b}]) \\
    &&+(-1)^{\lvert \mathbf{a'} \rvert+1}x(s^{\suspindex-1}[\mathbf{a'},\partial\mathbf{a}]\otimes[0,\mathbf{b}])\\
    &&-(-1)^{\lvert \mathbf{a'} \rvert+\lvert \mathbf{a} \rvert+1}x(s^{\suspindex-1}[\mathbf{a'},\mathbf{a}]\otimes[\mathbf{b}])\\
    &&+(-1)^{\lvert \mathbf{a'} \rvert\lvert \mathbf{a} \rvert+1}x(s^{\suspindex-1}[\mathbf{a'},\mathbf{a}]\otimes[0,\partial^{\circ}\mathbf{b}]).
    \end{array}
    \]
    The first two summands as well as the last one vanish again by \ref{CondAcross}, so the third summand also needs to vanish. This yields the desired equality. 

For $\lvert \mathbf{a}\rvert\geq 1, \lvert \mathbf{a'}\rvert\geq 0$ and $\lvert \mathbf{b}\rvert=1$,
    on the one hand
      \[
    \begin{array}{llll}
     \widetilde x\partial ([\mathbf{a'},\suspindex,\mathbf{a}] \otimes [b_0,b_1])   &=& \widetilde x([\partial\mathbf{a'},\suspindex,\mathbf{a}] \otimes [b_0,b_1])\\
     &&+(-1)^{\lvert \mathbf{a'}\vert+1}\widetilde x([\mathbf{a'},\mathbf{a}]\otimes [b_0,b_1])\\
     &&+(-1)^{\lvert \mathbf{a'}\vert+2}\widetilde x([\mathbf{a'},\suspindex,\partial\mathbf{a}] \otimes [b_0,b_1])\\
     &&+(-1)^{\lvert \mathbf{a'}\rvert+2+\lvert \mathbf{a}\rvert}\widetilde x([\mathbf{a'},\suspindex,\mathbf{a}] \otimes [b_0-b_1])&\\
         & =&(-1)^{\lvert \mathbf{a'}\vert+1}x(s^{\suspindex-1}[\mathbf{a'},\mathbf{a}] \otimes [b_0,b_1])\\
         &&+x(s^{\suspindex-1}[\mathbf{a'},\suspindex,\mathbf{a}] \otimes [0])&\\
         &&-x(s^{\suspindex-1}[\mathbf{a'},\suspindex,\mathbf{a}] \otimes [0]);&\\
    \end{array}
    \]
      on the other hand, we have
      \[
     \begin{array}{llll}
     \partial \widetilde x([\mathbf{a'},\suspindex,\mathbf{a}] \otimes [b_0,b_1])    &=& 0.&
     \end{array}
     \]
    Observe that from \ref{CondAcross}, we obtain
        \[
    \begin{array}{llll}
    0&=&\partial x(s^{\suspindex-1}[\mathbf{a'},\mathbf{a}]\otimes[0,b_0,b_1]) \\&=& x(s^{\suspindex-1}[\partial\mathbf{a'},\mathbf{a}]\otimes[0,b_0,b_1])\\
    &&+(-1)^{\lvert \mathbf{a'} \rvert+1}x(s^{\suspindex-1}[\mathbf{a'},\partial\mathbf{a}]\otimes[0,b_0,b_1])\\
    &&-(-1)^{\lvert \mathbf{a'} \rvert+\lvert \mathbf{a} \rvert+1}x(s^{\suspindex-1}[\mathbf{a'},\mathbf{a}]\otimes[b_0,b_1])\\
    &&+(-1)^{\lvert \mathbf{a'} \rvert+\lvert \mathbf{a} \rvert+1}x(s^{\suspindex-1}[\mathbf{a'},\mathbf{a}]\otimes[0,b_0])\\
    &&-(-1)^{\lvert \mathbf{a'} \rvert+\lvert \mathbf{a} \rvert+1}x(s^{\suspindex-1}[\mathbf{a'},\mathbf{a}]\otimes[0,b_1]).&
    \end{array}
    \]
    The first two summands as well as the last two vanish again by \ref{CondAcross}, so the third summand also needs to vanish. This yields the desired equality.

For $\lvert \mathbf{a}\rvert=0, \lvert \mathbf{a'}\rvert\geq 0$ and $\lvert \mathbf{b}\rvert\geq 2$; on the one hand we have
      \[
    \begin{array}{llll}
     \widetilde x(\partial ([\mathbf{a'},\suspindex,a] \otimes [\mathbf{b}]) )   &=& \widetilde x([\partial\mathbf{a'},\suspindex,a] \otimes [\mathbf{b}])\\
     &&+(-1)^{\lvert \mathbf{a'}\vert+1}\widetilde x([\mathbf{a'},a]\otimes [\mathbf{b}])\\
     &&+(-1)^{\lvert \mathbf{a'}\vert+2}\widetilde x([\mathbf{a'},\suspindex] \otimes [\mathbf{b}])\\
     &&+(-1)^{\lvert \mathbf{a'}\rvert+2}\widetilde x([\mathbf{a'},\suspindex,a] \otimes \partial^{\circ}[\mathbf{b}])&\\
         & =&(-1)^{\lvert \mathbf{a'}\vert+1}x(s^{\suspindex-1}[\mathbf{a'},a] \otimes [\mathbf{b}])&\\
         & &+(-1)^{\lvert \mathbf{a'}\vert+2}x([\mathbf{a'}] \otimes [0,\mathbf{b}]);
    \end{array}
    \]
      on the other hand we have 
      \[
     \begin{array}{llll}
     \partial \widetilde x([\mathbf{a'},\suspindex,a] \otimes [\mathbf{b}])    &=& 0.
     \end{array}
     \]
     Observe that by \ref{CondAcross}, we have
      \[
    \begin{array}{llll}
    0&=&\partial x(s^{\suspindex-1}[\mathbf{a'},a]\otimes[0,\mathbf{b}]) \\
    &=& x(s^{\suspindex-1}[\partial\mathbf{a'},a]\otimes[0,\mathbf{b}]) \\
    &&+ (-1)^{\lvert \mathbf{a'}\rvert+1}x(s^{\suspindex-1}[\mathbf{a'}]\otimes[0,\mathbf{b}])&\\
    &&-(-1)^{\lvert \mathbf{a'} \rvert+1}x(s^{\suspindex-1}[\mathbf{a'},a]\otimes[\mathbf{b}])\\
    &&-(-1)^{\lvert \mathbf{a'} \rvert+1}x(s^{\suspindex-1}[\mathbf{a'},a]\otimes[0,\partial^{\circ}\mathbf{b}]).&
    \end{array}
    \]
    The first and the last summands vanish using \ref{CondAcross} once again. This yields the desired equality.
     
For $\lvert \mathbf{a}\rvert=0, \lvert \mathbf{a'}\rvert\geq 0$, $\lvert \mathbf{b}\rvert=1$, on the one hand we have
      \[
    \begin{array}{llll}
     \widetilde x(\partial ([\mathbf{a'},\suspindex,a] \otimes [b_0,b_1]) )   &=& \widetilde x([\partial\mathbf{a'},\suspindex,a] \otimes [b_0,b_1])\\
     &&+(-1)^{\lvert \mathbf{a'}\vert+1}\widetilde x([\mathbf{a'},a]\otimes [b_0,b_1])\\
     &&+(-1)^{\lvert \mathbf{a'}\vert+2}\widetilde x([\mathbf{a'},\suspindex] \otimes [b_0,b_1])\\
     &&+(-1)^{\lvert \mathbf{a'}\rvert+2}\widetilde x([\mathbf{a'},\suspindex,a] \otimes [b_0-b_1])&\\
         & =&(-1)^{\lvert \mathbf{a'}\vert+1}x(s^{\suspindex-1}[\mathbf{a'},a] \otimes [b_0,b_1])&\\
         & &+(-1)^{\lvert \mathbf{a'}\vert+2}x([\mathbf{a'}] \otimes [0,b_0,b_1]) \\
     &&+x(s^{\suspindex-1}[\mathbf{a'},\suspindex,a] \otimes [0])\\
     &&-x(s^{\suspindex-1}[\mathbf{a'},\suspindex,a] \otimes [0]);
    \end{array}
    \]
      on the other hand we have 
      \[
     \begin{array}{llll}
     \partial \widetilde x([\mathbf{a'},\suspindex,a] \otimes [b_0,b_1])    &=& 0.
     \end{array}
     \]
     Observe that by \ref{CondAcross}, we have
      \[
    \begin{array}{llll}
    0&=&\partial x(s^{\suspindex-1}[\mathbf{a'},a]\otimes[0,b_0,b_1])\\
    &=& x(s^{\suspindex-1}[\partial\mathbf{a'},a]\otimes[0,b_0,b_1]) \\
    &&+ (-1)^{\lvert \mathbf{a'}\rvert+1}x(s^{\suspindex-1}[\mathbf{a'}]\otimes[0,b_0,b_1])&\\
    &&-(-1)^{\lvert \mathbf{a'} \rvert+1}x(s^{\suspindex-1}[\mathbf{a'},a]\otimes[b_0,b_1])\\
    &&-(-1)^{\lvert \mathbf{a'} \rvert+1}x(s^{\suspindex-1}[\mathbf{a'},a]\otimes[0,b_0])&\\
    &&+(-1)^{\lvert \mathbf{a'} \rvert+1}x(s^{\suspindex-1}[\mathbf{a'},a]\otimes[0,b_1]).
    \end{array}
    \]
    The first and the last two summands vanish using \ref{CondAcross} once again. This yields the desired equality.
     
    
For $\lvert \mathbf{a}\rvert\geq 1, \lvert \mathbf{a'}\rvert=-1$, $\lvert \mathbf{b} \rvert\geq 2$, on the one hand we have
      \[
    \begin{array}{llll}
     \widetilde x(\partial ([\suspindex,\mathbf{a}] \otimes [\mathbf{b}]) )   &=& \widetilde x([\mathbf{a}] \otimes [\mathbf{b}])&\\
     &&-\widetilde x([\suspindex,\partial\mathbf{a}] \otimes [\mathbf{b}])\\
     &&+(-1)^{1+\lvert \mathbf{a}\rvert}\widetilde x([\suspindex,\mathbf{a}] \otimes \partial^{\circ}[\mathbf{b}])&\\
         & =&x(s^{\suspindex-1}[\mathbf{a}] \otimes [\mathbf{b}]);&\\
    \end{array}
    \]
      on the other hand we have 
      \[
     \begin{array}{llll}
     \partial \widetilde x([\suspindex,\mathbf{a}] \otimes [\mathbf{b}])    &=& 0.
         \\
     \end{array}
     \]
     Observe that by \ref{CondBelow}, the two results coincide. 

   
For $\lvert \mathbf{a}\rvert\geq 1, \lvert \mathbf{a'}\rvert=-1$, $\lvert \mathbf{b} \rvert= 1$, on the one hand we have
      \[
    \begin{array}{llll}
     \widetilde x(\partial ([\suspindex,\mathbf{a}] \otimes [b_0,b_1]) )   &=& \widetilde x([\mathbf{a}] \otimes [b_0,b_1])\\
     &&-\widetilde x([\suspindex,\partial\mathbf{a}] \otimes [b_0,b_1])\\
     &&+(-1)^{1+\lvert \mathbf{a}\rvert}\widetilde x( ([\suspindex,\mathbf{a}] \otimes [b_0]) )&\\
     &&-(-1)^{1+\lvert \mathbf{a}\rvert}\widetilde x(\partial ([\suspindex,\mathbf{a}] \otimes [b_1]) )&\\
         & =&x([s^{\suspindex-1}(\mathbf{a})] \otimes [b_0,b_1])\\
         &&
         +(-1)^{1+\lvert \mathbf{a}\rvert}x(s^{\suspindex-1}[\suspindex,\mathbf{a}] \otimes [0])\\
         &&
         -(-1)^{1+\lvert \mathbf{a}\rvert}x(s^{\suspindex-1}[\suspindex,\mathbf{a}] \otimes [0]);
    \end{array}
    \]
      on the other hand we have 
      \[
     \begin{array}{llll}
     \partial \widetilde x([\suspindex,\mathbf{a}] \otimes [b_0,b_1])    &=& 0.
         \\
     \end{array}
     \]
     Then \ref{CondBelow} yields the desired equality. 
For $\lvert \mathbf{a}\rvert=0$, $\lvert \mathbf{a'}\rvert=-1$, $\lvert \mathbf{b} \rvert\geq 2$, on the one hand we have
      \[
    \begin{array}{llll}
     \widetilde x(\partial ([\suspindex,a] \otimes [\mathbf{b}]) )   &=& \widetilde x([a] \otimes [\mathbf{b}])-\widetilde x([\suspindex]\otimes [\mathbf{b}])\\
     &&+(-1)^{1}\widetilde x([\suspindex,a] \otimes \partial^{\circ}[\mathbf{b}])&\\
         & =&x(s^{\suspindex-1}[a] \otimes [\mathbf{b}]);&\\
         & &&
    \end{array}
    \]
      on the other hand we have
      \[
     \begin{array}{llll}
     \partial \widetilde x([\suspindex,a] \otimes [\mathbf{b}])    &=&0.
         \\
     \end{array}
     \]
     Using \ref{CondBelow} once again, we obtain the desired equality.

For $\lvert \mathbf{a}\rvert=0$, $\lvert \mathbf{a'}\rvert=-1$, $\lvert \mathbf{b} \rvert=1$;
 on the one hand
      \[
    \begin{array}{llll}
     \widetilde x(\partial ([\suspindex,a] \otimes [b_0,b_1]) )   &=& \widetilde x([a] \otimes [b_0,b_1])-\widetilde x([r]\otimes [b_0,b_1])\\
     &&+\widetilde x([\suspindex,a] \otimes [b_1])-\widetilde x([r,a]\otimes [b_0])\\
         & =&x([s^{\suspindex-1}(a)] \otimes [b_0,b_1])\\
         &&+x(s^{\suspindex-1}[\suspindex,a] \otimes [0])\\
         &&-x(s^{\suspindex-1}[\suspindex,a] \otimes [0]);&\\
    \end{array}
    \]
      on the other hand we have 
      \[
     \begin{array}{llll}
     \partial \widetilde x([\suspindex,a] \otimes [b_0,b_1])    &=& 0.
         \\
     \end{array}
     \]
     Using \ref{CondBelow} once again yields the desired equality. 

     \item In this case, all constituents can be seen to be $0$ in a straightforward manner.
\end{enumerate}
Since we treated all possible cases, we conclude that $\widetilde x$ is indeed a chain map. 
\end{proof}

\bibliographystyle{amsalpha}
\bibliography{ref}

\end{document}